\documentclass[twoside,12pt]{article}
\usepackage{latexsym}
\usepackage{amsmath}
\usepackage{amssymb}
\usepackage{cases}
\usepackage{ascmac}
\usepackage{mathrsfs}
\usepackage{amsthm}
\usepackage{cite}
\usepackage[usenames,dvipsnames]{color}

\title{Time periodic solutions of {C}ahn--{H}illiard system with dynamic boundary conditions }

\author{%Takeshi Fukao\\
Taishi Motoda\\
Graduate School of Education, 
Kyoto University of Education\\
1~Fujinomori, Fukakusa, Fushimi-ku, Kyoto~612-8522 Japan\\
E-mail: \texttt{motoda.math@gmail.com}}
 
\date{}

\newcommand\testopari{\sc Taishi Motoda}
\newcommand\testodispari{\sc Time periodic solutions of {C}ahn--{H}illiard system }
\begin{document}

\maketitle

\begin{abstract}
The existence problem for {C}ahn--{H}illiard system with dynamic boundary conditions and time periodic conditions is discussed. 
We apply the abstract theory of evolution equations by using viscosity approach 
and the Schauder fixed point theorem in the level of approximate ploblem. 
One of the key point is the assumption for maximal monotone graphs with respect to their domains. 
Thanks to this, we obtain the existence result of the weak solution by using the passage to the limit. 

\vspace{2mm}

\vspace{2mm}

\end{abstract}

%%%%% Section 1. %%%%%
\section{Introduction}
\setcounter{equation}{0}
\indent

In this paper, we consider the following {C}ahn--{H}illiard system with dynamic boundary condition and time periodic condition, 
say {\rm (P)}, which consists of the following: 
\begin{eqnarray}
&&\displaystyle \frac{\partial u}{\partial t}-\Delta \mu =0 \quad {\rm in~}Q:=\Omega \times (0, T), \label{eq1} \\
&&\mu =-\kappa _{1}\Delta u+\xi +\pi (u)-f, \quad \xi \in \beta (u) \quad {\rm in~}Q, \label{eq2} \\
&&u_{\Gamma }=u_{|_{\Gamma }}, \quad \mu _{\Gamma }=\mu _{|_{\Gamma }} \quad {\rm on~}\Sigma :=\Gamma \times (0, T), \label{eq3} \\
&&\displaystyle \frac{\partial u_{\Gamma }}{\partial t}+\partial _{\boldsymbol{\nu }}\mu -\Delta _{\Gamma }\mu _{\Gamma }=0 \quad {\rm on~}\Sigma , \label{eq4} \\
&&\mu _{\Gamma }=\kappa _{1}\partial _{\boldsymbol{\nu }}u-\kappa _{2}\Delta _{\Gamma }u_{\Gamma }+\xi _{\Gamma }+\pi _{\Gamma }(u_{\Gamma })-f_{\Gamma }, \quad \xi _{\Gamma }\in \beta _{\Gamma }(u_{\Gamma }) \quad {\rm on~}\Sigma , \label{eq5} \\
&&u(0)=u(T) \quad {\rm in~}\Omega , \quad u_{\Gamma }(0)=u_{\Gamma }(T) \quad {\rm on~}\Gamma \label{eq6}
\end{eqnarray}
where $0< T< +\infty$, $\Omega $ is a bounded domain of $\mathbb{R}^{d}$ $(d=2,3)$ 
with smooth boundary $\Gamma :=\partial \Omega $, 
$\kappa _{1}, \kappa _{2}$ are positive constants, 
$\partial _{\boldsymbol{\nu}}$ is the outward normal derivative on $\Gamma $, 
$u _{|_{\Gamma }}, \mu _{|_{\Gamma }}$ stand for the trace of $u$ and $\mu $ to $\Gamma $, respectively, 
$\Delta $ is the {L}aplacian, 
$\Delta _{\Gamma }$ is the {L}aplace--{B}eltrami operator 
(see, e.g., \cite{Gri09}), and
$f: Q\to \mathbb{R}$, 
$f_{\Gamma }: \Sigma \to \mathbb{R}$ are given data. 
Moreover, in the nonlinear diffusion terms, 
$\beta , \beta _{\Gamma }: \mathbb{R}\to 2^{\mathbb{R}}$ are maximal monotone graphs and 
$\pi , \pi _{\Gamma }: \mathbb{R}\to \mathbb{R}$ are {L}ipschitz perturbations. 

The {C}ahn--{H}illiard equation \cite{CH58} is a description of 
mathematical model for phase separation, e.g., 
the phenomenon of separating into two phases from homogeneous composition, 
the so-called spinodal decomposition. 
In (\ref{eq1})--(\ref{eq2}), $u$ is the order parameter and 
$\mu $ is the chemical potential. 
Moreover, it is well known that the {C}ahn--{H}illiard equation 
is characterized by the nonlinear term $\beta +\pi $. 
It play an important role as derivative of double-well potential $W$. 
As pioneering results, {K}enmochi, {N}iezg\'odka and {P}aw\l ow study 
the {C}ahn--{H}illiard equation with constraint by subdifferential operator approach \cite{DK99}. 
In addition, {K}ubo investigates the {C}ahn--{H}illiard equation \cite{Kub12}. 
Recently, in \cite{Kaj17}, it is discussed the strong solution of the {C}ahn--{H}illiard equation in bounded domains with permeable
and non-permeable walls in maximal $L_{p}$ regularity spaces. 

In terms of (\ref{eq3})--(\ref{eq5}), 
we consider the dynamic boundary condition as being $u_{\Gamma }, \mu _{\Gamma }$ unknown functions
such that satisfied trace condition (\ref{eq3}). 
The dynamic boundary condition was treated in recent years, 
for example, for the Stefan proplem \cite{Aik93, Aik95, Aik96, Fuk16-1}, 
wider the degenerate parabolic equation \cite{Fuk16-2, FM18} 
and the {C}ahn--{H}illiard equation \cite{CF15, CGS14, GMS09, GMS13, LW17}. 
To the best our knowledge, the type of dynamic boundary condition on the {C}ahn--{H}illiard equation like (P) is formulated in \cite{GMS11}. 
As you can see, we consider the same type of equations (\ref{eq1})--(\ref{eq2}) on the boundary. 
In other words, (\ref{eq1})--(\ref{eq5}) is a transmission problem connecting $\Omega $ and $\Gamma $.

The nonlinear term $\beta _{\Gamma }+\pi _{\Gamma }$ 
is also the derivative of double-well potential $W_{\Gamma }$. 
As a well-known example, $W=W_{\Gamma }=(1/4)(r^{2}-1)^{2}$, namely, $W'=W_{\Gamma }'=r^{3}-r$ for $r\in \mathbb{R}$. 
This is called the prototype double well potential. 
Note that we do not have to take differnt nonlinear terms $W'$ and $W_{\Gamma }'$ in $\Omega $ and on $\Gamma $, respectively. 
However, in problem (P), we treat it differently to generalize. 
In this case, it is necessary to assume the compatibility condition (see e.g., \cite{CC13, CF15}), stated (A4). 
The other example of $W, W_{\Gamma }$ is stated later.

About the problem (P), {C}olli and {F}ukao study the 
{C}ahn--{H}illiard system with daynamic boundary condition and initial value condition \cite{CF15}. 
They set a function space that the total mass is zero. 
This idea is arised from the property of dynamic boundary condition. 
The property is called the toal mass conservation.

Moreover, focusing on (\ref{eq6}), 
the study respect to existence of time periodic solutions of the {C}ahn--{H}illiard equation is not much. 
For example, \cite{YJ11, LYRJ08, LWZ16, WZ16}. 
In particular, 
{W}ang and {Z}heng discuss the existence of 
time periodic solution of the {C}ahn--{H}illiard equation with {N}eumann boundary condition \cite{WZ16}. 
They employ the method of \cite{AS11}. 
%, especialy, 
%use the viscosity approach, the {S}chauder fixed point theorem in the level of approximate problem 
%and the theory of subdifferential operator (see e.g., \cite{Kem81}). 
Note that they impose two assumptions for a maximal monotone graph, 
specifically, restricted domains and the following growth condition for maximal monotone graph $\beta $ : 
\begin{equation*}
\widehat{\beta }(r)\geq cr^{2} \quad {\rm for~all~}r\in \mathbb{R}, 
\end{equation*}
for some positive constant $c$. 
However, the above assumption is too restrictive for some physical applications.

%Moreover, in \cite{LWZ16}, 
%the existence of time periodic solution for the {C}ahn--{H}illiard equation with {N}eumann boundary condition is discussed. 

In this paper, following the method of \cite{WZ16}, 
we apply the abstract theory of evolution equations by using the viscosity approach 
and the {S}chauder fixed point theorem in the level of approximate problem. 
Moreover, by virtue of the viscosity approach, 
we can apply the abstract result \cite{AS11}. 
Noting that, in \cite{CF15}, they do not impose the assumption of growth condition 
for maximal monotone graphs $\beta $ and $\beta _{\Gamma }$. 
%The main objective of this paper is to improve 
Therefore, by setting the appropriate convex functional and using the {P}oincar\'e--{W}irtinger inequality, 
we can also avoid imposing the growth condition 
even though we consider the time periodic problem. 
Thanks to this result, we can choose wider kinds of nonlinear diffusion terms $\beta +\pi $ and $\beta _{\Gamma }+\pi _{\Gamma }$. 
However, respect to the assumption of restricted domains (see (A5)), this is a essential point to slove the problem (P). 
Nevertheless, we can infer that the assumption (A5) is not necessary 
when we choose the prototype double well potential. 
We sketch it Remark 4.1.

The present paper proceeds as follows.

In Section 2, main theorem and definition of solution are stated. 
At first, we prepare the notation used in this paper and set appropriate 
function spaces. 
Next, we introduce the definition of weak solution of {\rm (P)} %and {\rm (P)}$_{\varepsilon }$
and the main theorems are given there. 
Also, we give the example of double-well potentials. 

In Section 3, in order to prove convergence theorem, at first, 
we set convex functionals and consider approximate problems. 
Next, we obtain the solution of (P)$_{\varepsilon }$ by using the {S}chauder fixed point theorem. 
Finally, we deduce uniform estimates of the solution of (P)$_{\varepsilon }$.

In Section 4, we prove the existence of weak solutions by passing to the limit $\varepsilon \to 0$. 
%We also discuss the uniqueness of solutions. 
Finally, we sketch the case that we choose the derivative of prototype double well potential as nonlinear diffusion terms. 

%In Section 5, we prove the error estimate. 

A detailed index of sections and subsections follows.

\begin{itemize}
 \item[1.] Introduction
 \item[2.] Main results
\begin{itemize}
 \item[2.1.] Notation
 \item[2.2.] Definition of the solution and main theorem
\end{itemize}
 \item[3.] Approximate problem and uniform estimates
\begin{itemize}
 \item[3.1.] Abstract formulation
 \item[3.2.] Approximate problem for (P)
 \item[3.3.] Uniform estimates
\end{itemize}
 \item[4.] Proof of convergence theorem
%\begin{itemize}
% \item [4.1.] Passage to the limit $\lambda \searrow 0$
% \item [4.2.] Passage to the limit $\varepsilon \searrow 0$
%\end{itemize}
\end{itemize}

%%%%% Section 2. %%%%%
\section{Main results}
\setcounter{equation}{0}

%%%%% 2.1 %%%%%
\subsection{Notation}
\indent

We introduce the spaces 
$H:=L^{2}(\Omega )$, 
$H_{\Gamma }:=L^{2}(\Gamma )$, 
$V:=H^{1}(\Omega )$, 
$V_{\Gamma }:=H^{1}(\Gamma )$
with standard norms 
$|\cdot |_{H}$, 
$|\cdot |_{H_{\Gamma }}$, 
$|\cdot |_{V}$, 
$|\cdot |_{V_{\Gamma }}$
and inner products 
$(\cdot , \cdot )_{H}$, 
$(\cdot , \cdot )_{H_{\Gamma }}$, 
$(\cdot , \cdot )_{V}$, 
$(\cdot , \cdot )_{V_{\Gamma }}$, respectively. 
Moreover, we set $\boldsymbol{H}:=H\times H_{\Gamma }$ and
\begin{equation*}
	\boldsymbol{V}:= 
	\bigl\{ 
	\boldsymbol{z}:=(z, z_{\Gamma }) 
	\in V \times V_{\Gamma } \ : \ z_{|_\Gamma }=z_{\Gamma } \ {\rm a.e.~on}~\Gamma 
	\bigr\}. 
\end{equation*}
$\boldsymbol{H}$ and $\boldsymbol{V}$ are then {H}ilbert spaces with inner products
\begin{gather*}
	(\boldsymbol{u}, \boldsymbol{z})_{\boldsymbol{H}}
	:=(u, z)_{H}+(u_{\Gamma }, z_{\Gamma })_{H_{\Gamma }} 
	\quad {\rm for~all~} 
	\boldsymbol{u}:=(u, u_{\Gamma }), \boldsymbol{z}:=(z, z_{\Gamma })\in \boldsymbol{H}, 
	\\
	(\boldsymbol{u}, \boldsymbol{z})_{\boldsymbol{V}}
	:=(u, z)_{V}+(u_{\Gamma }, z_{\Gamma })_{V_{\Gamma }} 
	\quad \quad {\rm for~all~} \boldsymbol{u}
	:=(u, u_{\Gamma }), \boldsymbol{z}:=(z, z_{\Gamma })\in \boldsymbol{V}.
\end{gather*}
Note that 
$\boldsymbol{z} \in \boldsymbol{V}$ 
implies that the second component $z_{\Gamma }$ of $\boldsymbol{z}$ 
is equal to the trace of the first component $z$ of $\boldsymbol{z}$ on $\Gamma $, and 
$\boldsymbol{z} \in \boldsymbol{H}$ implies that 
$z\in H$ and $z_{\Gamma }\in H_{\Gamma }$ are independent. 
Throughout this paper,
we use the bold letter $\boldsymbol{u}$ to represent 
the pair corresponding to the letter; i.e., $\boldsymbol{u}:=(u, u_{\Gamma })$.

Let $m:\boldsymbol{H}\to \mathbb{R}$ be the mean function defined by 
\begin{equation*}
	m(\boldsymbol{z})
	:=\frac{1}{|\Omega |+|\Gamma |}
	\left\{\int_{\Omega }zdx+\int_{\Gamma }z_{\Gamma }d\Gamma \right\}
	\quad {\rm for~all~} \boldsymbol{z} \in \boldsymbol{H},
\end{equation*}
where $|\Omega |:=\int_{\Omega }1dx, |\Gamma |:=\int_{\Gamma }1d\Gamma $. 
Then, we define 
$\boldsymbol{H}_{0}:=\{\boldsymbol{z}\in \boldsymbol{H};m(\boldsymbol{z})=0\}$, $\boldsymbol{V}_{0}:=\boldsymbol{V}\cap \boldsymbol{H}_{0}$. 
Moreover, 
$\boldsymbol{V}^{*}, \boldsymbol{V}_{0}^{*}$ denote the dual spaces of $\boldsymbol{V}, \boldsymbol{V}_{0}$, respectively; 
the duality pairing between $\boldsymbol{V}_{0}^{*}$ and $\boldsymbol{V}_{0}$ is denoted 
$\langle  \cdot , \cdot \rangle _{\boldsymbol{V}_{0}^{*}, \boldsymbol{V}_{0}}$. 
We define the norm of $\boldsymbol{H}_{0}$
by $|\boldsymbol{z}|_{\boldsymbol{H}_{0}}:=|\boldsymbol{z}|_{\boldsymbol{H}}$ 
for all $\boldsymbol{z}\in \boldsymbol{H}_{0}$. 
Now, we define the 
bilinear form $a(\cdot , \cdot ):\boldsymbol{V}\times \boldsymbol{V}\to \mathbb{R}$ by
\begin{equation*} 
	a(\boldsymbol{u}, \boldsymbol{z})
	:= \kappa _{1}\int_{\Omega } \nabla u\cdot \nabla zdx
	+ \kappa _{2}\int_{\Gamma } \nabla _{\Gamma } u_{\Gamma }\cdot \nabla _{\Gamma }z_{\Gamma }d\Gamma 
	\quad {\rm for~all~} \boldsymbol{u}, \boldsymbol{z}\in \boldsymbol{V}. 
\end{equation*}
Then, for all $\boldsymbol{z}\in \boldsymbol{V}_{0}$, 
$|\boldsymbol{z}|_{\boldsymbol{V}_{0}}
:=\sqrt{a(\boldsymbol{z}, \boldsymbol{z})}$ is the norm of $\boldsymbol{V}_0$. 
Also, for all $\boldsymbol{z}\in \boldsymbol{V}_{0}$, 
we let $\boldsymbol{F}:\boldsymbol{V}_{0}\to \boldsymbol{V}_{0}^{*}$ be 
the duality mapping defined by
\begin{equation*}
	\langle \boldsymbol{F}\boldsymbol{z}, 
	\tilde{\boldsymbol{z}}
	\rangle _{\boldsymbol{V}_{0}^{*}, \boldsymbol{V}_{0}}
	:=a(\boldsymbol{z}, \tilde{\boldsymbol{z}}) 
	\quad {\rm for~all~}
	\tilde{\boldsymbol{z}}\in \boldsymbol{V}_{0}.
\end{equation*}
Then, the following the {P}oincar\'e--{W}irtinger inequality holds: there exists a positive constant $c_{\rm P}$ such that 
\begin{equation}
	\label{PW0}
	 |\boldsymbol{z}|_{\boldsymbol{V}}^{2} 
	\le c_{\rm P} 
	 |\boldsymbol{z}|_{\boldsymbol{V}_{0}}^{2}
	\quad {\rm for~all~} 
	\boldsymbol{z}\in \boldsymbol{V}
\end{equation}
(see \cite[Lemma C]{CF15}). 
Moreover, we define the inner product of $\boldsymbol{V}_{0}^{*}$ by
\begin{equation*}
	(\boldsymbol{z}^{*}, \tilde{\boldsymbol{z}}^{*})_{\boldsymbol{V}_{0}^{*}}
	:=\langle  \boldsymbol{z}^{*}, \boldsymbol{F}^{-1}\tilde{\boldsymbol{z}}^{*}
	\rangle _{\boldsymbol{V}_{0}^{*}, \boldsymbol{V}_{0}} 
	\quad 
	{\rm for~all~} 
	\boldsymbol{z}^{*}, \tilde{\boldsymbol{z}}^{*}\in \boldsymbol{V}_{0}^{*}.
\end{equation*}
Also, we define the projection $\boldsymbol{P}: \boldsymbol{H}\to \boldsymbol{H}_{0}$ by 
\begin{equation*}
\boldsymbol{P}\boldsymbol{z}:=\boldsymbol{z}-m(\boldsymbol{z})\boldsymbol{1}%=\bigl(z-m(\boldsymbol{z}), z_{\Gamma }-m(\boldsymbol{z})\bigr)
\quad {\rm for~all~}\boldsymbol{z}\in \boldsymbol{H}, 
\end{equation*}
where $\boldsymbol{1}:=(1, 1)$. 
Then, since $\boldsymbol{P}$ is a linear bounded operator, note that the following property holds: 
let $\{\boldsymbol{z}_{n}\}_{n\in \mathbb{N}}$ be a sequence in $\boldsymbol{H}$ such that 
$\boldsymbol{z}_{n}\to \boldsymbol{z}$ weakly in $\boldsymbol{H}$ for some $\boldsymbol{z}$, 
then we infer that 
\begin{equation}\label{projection}
\boldsymbol{P}\boldsymbol{z}_{n}\to \boldsymbol{P}\boldsymbol{z} \quad {\rm weakly~in~}\boldsymbol{H}_{0}, 
\end{equation}
because $\boldsymbol{P}$ is a single-valued operator. 
Then, we have
$\boldsymbol{V}_{0}\hookrightarrow \hookrightarrow \boldsymbol{H}_{0}\hookrightarrow \hookrightarrow \boldsymbol{V}_{0}^{*}$, where 
``$\hookrightarrow \hookrightarrow $'' stands for 
compact embedding (see \cite[Lemmas A and B]{CF15}).

%%%%% 2.2 %%%%%
\subsection{Definition of the solution and main theorem}
\indent

In this subsection we define our solution for {\rm (P)} and then 
we state the main theorem. 

Firstly, from (\ref{eq1}), (\ref{eq4}) and the property of dynamic boundary condition, it holds the following total mass conservation: 
\begin{equation*}
\int_{\Omega }u(t)dx+\int_{\Gamma }u_{\Gamma }(t)d\Gamma =\int_{\Omega }u(0)dx+\int_{\Gamma }u_{\Gamma }(0)d\Gamma :=m_{0} \quad {\rm for~all~}t\in [0, T]. 
\end{equation*}

To define solutions, we use the following notation: 
the variable $\boldsymbol{v}:=\boldsymbol{u}-m_{0}\boldsymbol{1}$; 
the datum $\boldsymbol{f}:=(f, f_{\Gamma })$; 
the function $\boldsymbol{\pi }(\boldsymbol{z}):=(\pi (z), \pi _{\Gamma }(z_{\Gamma }))$ for $\boldsymbol{z}\in \boldsymbol{H}$. 
Moreover, we set the space $\boldsymbol{W}:=H^{2}(\Omega )\times H^{2}(\Gamma )$.

\paragraph{Definition 2.1.}
{\it The triplet 
$(\boldsymbol{v}, \boldsymbol{\mu }, \boldsymbol{\xi })$ 
is called the weak solution of {\rm (P)} if 
\begin{gather*}
	\boldsymbol{v}\in H^{1}(0, T;\boldsymbol{V}_{0}^{*})\cap L^{\infty }(0, T;\boldsymbol{V}_{0})\cap L^{2}(0, T;\boldsymbol{W}), \\
	\boldsymbol{\mu }\in L^{2}(0, T;\boldsymbol{V}), \\
	\boldsymbol{\xi }=(\xi , \xi _{\Gamma })\in L^{2}(0, T;\boldsymbol{H}), 
\end{gather*}
and they satisfy
\begin{gather}
\bigl\langle \boldsymbol{v}'(t), \boldsymbol{z}\bigr\rangle _{\boldsymbol{V}_{0}^{*}, \boldsymbol{V}_{0}}+a\bigl(\boldsymbol{\mu }(t), \boldsymbol{z}\bigr)=0 \quad {\rm for~all~}z\in \boldsymbol{V}_{0} \label{weak1} \\
\bigl(\boldsymbol{\mu }(t), \boldsymbol{z}\bigr)_{\boldsymbol{H}}=a\bigl(\boldsymbol{v}(t), \boldsymbol{z}\bigr)
+\bigl(\boldsymbol{\xi }(t)-m\bigl(\boldsymbol{\xi }(t)\bigr)+\boldsymbol{\pi }\bigl(\boldsymbol{v}(t)+m_{0}\boldsymbol{1}\bigr)-\boldsymbol{f}, \boldsymbol{z}\bigr)_{\boldsymbol{H}} \quad {\rm for~all~}z\in \boldsymbol{V} \label{weak2}
\end{gather}
for a.a.\ $t\in (0, T)$, and
\begin{gather*}
\xi \in \beta (v+m_{0}) \quad {\it a.e.~in~}Q, \quad \xi _{\Gamma }\in \beta _{\Gamma }(v_{\Gamma }+m_{0}) \quad {\it a.e.~on~}\Sigma 
\end{gather*}
with}
\begin{equation}
	\boldsymbol{v}(0)=\boldsymbol{v}(T) \quad {\it in~}\boldsymbol{H}_{0}. \label{ic}
\end{equation}

\paragraph{Remark 2.1.}
We can see that $\boldsymbol{\mu }:=(\mu , \mu _{\Gamma })$ satisfies 
\begin{gather*}
\mu =-\kappa _{1}\Delta u+\xi -m\bigl(\boldsymbol{\xi }\bigr)+\pi (u)-f \quad {\rm a.e.~in~}Q, \\
\mu _{\Gamma }=\kappa _{1}\partial _{\boldsymbol{\nu }}u-\kappa _{2}\Delta _{\Gamma }u_{\Gamma }
+\xi _{\Gamma }-m\bigl(\boldsymbol{\xi }\bigr)+\pi _{\Gamma }(u_{\Gamma })-f_{\Gamma } \quad {\rm a.e.~on~}\Sigma , 
\end{gather*}
where $u=v+m_{0}$ and $u_{\Gamma }=v_{\Gamma }+m_{0}$, because of the regularity $\boldsymbol{v}\in L^{2}(0, T; \boldsymbol{W})$.

\paragraph{Remark 2.2.}
In (\ref{weak2}), this is different from the following definition of \cite[Difinition~2.1]{CF15}: 
\begin{gather}\label{cfdef}
\bigl(\boldsymbol{\mu }(t), \boldsymbol{z}\bigr)_{\boldsymbol{H}}=a\bigl(\boldsymbol{v}(t), \boldsymbol{z}\bigr)
+\bigl(\boldsymbol{\xi }(t)+\boldsymbol{\pi }\bigl(\boldsymbol{v}(t)+m_{0}\boldsymbol{1}\bigr)-\boldsymbol{f}, \boldsymbol{z}\bigr)_{\boldsymbol{H}} \quad {\rm for~all~}z\in \boldsymbol{V}
\end{gather}
for a.a.\ $t\in (0, T)$. 
However, by setting $\widetilde{\boldsymbol{\mu }}:=\boldsymbol{\mu }+m\bigl(\boldsymbol{\xi }\bigr)\boldsymbol{1}$, 
$\widetilde{\boldsymbol{\mu }}$ satisfies $\widetilde{\boldsymbol{\mu }}\in L^{2}(0, T; \boldsymbol{V})$ and (\ref{cfdef}). 
Hence, in other words, we can employ (\ref{cfdef}) as definition of (P) replaced by (\ref{weak2}). \\
%\paragraph{Definition 2.2.}
%{\it The triplet 
%$(\boldsymbol{v}_{\varepsilon }, \boldsymbol{\mu }_{\varepsilon }, \boldsymbol{\xi }_{\varepsilon })$ 
%is called the weak solution of {\rm (P)}}$_{\varepsilon }$ {\it if 
%\begin{gather}
%	u_{\varepsilon }\in H^{1}(0, T;V^{*})\cap L^{\infty }(0, T;V), \label{3.1}\\
%	\mu _{\varepsilon }\in L^{2}(0, T;V), \quad \xi _{\varepsilon }\in L^{2}(0, T;H), \label{3.2}\\
%	\xi _{\varepsilon }\in \beta (u_{\varepsilon }) \quad {\it a.e.~in~}Q \label{maximala}
%\end{gather}
%and they satisfy
%\begin{gather}
%\bigl\langle u_{\varepsilon }'(t), z\bigr\rangle _{V^{*}, V}+\int_{\Omega }\nabla \mu _{\varepsilon }(t)\cdot \nabla zdx+\kappa \int_{\Gamma }\mu _{\varepsilon }zd\Gamma =0 \quad {\it for~all~}z\in V \label{2.13}\\
%\mu _{\varepsilon }=-\varepsilon \Delta u_{\varepsilon }(t)+\xi _{\varepsilon }(t)+\pi _{\varepsilon }(u_{\varepsilon }(t))-f(t) \quad {\it in~}H \label{2.14}
%\end{gather}
%for a.a.\ $t\in (0, T)$, with}
%\begin{equation}
%	u(0)=u_{0} \quad {\it a.e.~in}~\Omega . \label{2.16}
%\end{equation}

%%%assumption%%%%%

We assume that
\begin{enumerate}
 \item[(A1)] $\boldsymbol{f}\in L^{2}(0, T; \boldsymbol{V})$ and $\boldsymbol{f}(t)=\boldsymbol{f}(t+T)$ for a.a. $t\in (0, T)$. 
 \item[(A2)] $\pi , \pi _{\Gamma }: \mathbb{R}\to \mathbb{R}$ 
are locally {L}ipschitz continuous functions with {L}ipschitz contants $L, L_{\Gamma }$, respectively; 
 \item[(A3)] $\beta , \beta _{\Gamma }:\mathbb{R}\to 2^{\mathbb{R}}$ are maximal monotone graghs, which is
the subdifferential 
$$\beta =\partial _{\mathbb{R}}\widehat{\beta }, \quad \beta _{\Gamma }=\partial _{\mathbb{R}}\widehat{\beta }_{\Gamma }$$
of some proper lower semicontinuous convex functions 
$\widehat{\beta }, \widehat{\beta }_{\Gamma }: \mathbb{R}\to [0, +\infty ]$ 
satisfying $\widehat{\beta }(0)=\widehat{\beta }_{\Gamma }(0)=0$ with domains $D(\beta )$ and $D(\beta _{\Gamma })$, respectively. 
 \item[(A4)] $D(\beta _{\Gamma })\subseteq D(\beta )$ and there exist positive constant $\rho , c_{0}>0$ such that 
\begin{equation}\label{ccbeta}
|\beta ^{\circ }(r)|\leq \rho |\beta _{\Gamma }^{\circ }(r)|+c_{0} \quad {\rm for~all~}r\in D(\beta _{\Gamma }): 
\end{equation}
 \item[(A5)]$D(\beta ), D(\beta _{\Gamma })$ are bounded domains with 
non-empty interior and , {\rm i.e.}, 
$\overline{D(\beta )}=[\sigma _{*}, \sigma ^{*}]$ and $\overline{D(\beta _{\Gamma })}=[\sigma _{\Gamma *}, \sigma _{\Gamma }^{*}]$
for some constants $\sigma _{*}, \sigma ^{*}, \sigma _{\Gamma *}$ and $\sigma _{\Gamma }^{*}$ 
with $-\infty <\sigma _{*}\leq \sigma _{\Gamma *}<\sigma _{\Gamma }^{*}\leq \sigma ^{*}<\infty . $
\end{enumerate}
In particular, (A3) yields $0\in \beta (0)$. 
The assumption (A5) is not imposed \cite{CF15}. 
However, it is essential to obtain uniform estimates in Section 3. 
This is the greatest difficulties of time periodic problem. 
Also, the assumption of compatibility of $\beta $ and $\beta _{\Gamma }$ (A4) is the same as in \cite{CC13, CF15}. 
Now, we give the example respect to nonlinear diffusion terms under the above assumptions: \\
~\\
$\bullet $ $\beta (r)=\beta_{\Gamma }(r)=(\alpha _{1}/2)\ln ((1+r)/(1-r))$, 
$\pi (r)=\pi _{\Gamma }(r)=-\alpha _{2}r$ for all $r\in D(\beta )=D(\beta _{\Gamma })=(-1. 1)$ 
and $0<\alpha _{1}<\alpha _{2}$. 
for the logarithmic double well potential 
$W(r)=W_{\Gamma }(r)=(\alpha _{1}/2)\{(1-r)\ln ((1-r)/2)+(1+r)\ln ((1+s)/2)\}+(\alpha _{2}/2)(1-r^{2})$. 
The condition $\alpha _{1}<\alpha _{2}$ ensures that $W, W_{\Gamma }$ have double-well forms 
(see e.g., \cite{CMZ11}). \\
$\bullet $ $\beta (r)=\beta _{\Gamma }(r)=\partial I_{[-1, 1]}(r)$, 
$\pi (r)=\pi _{\Gamma }(r)=-r$ for all $r\in D(\beta )=D(\beta _{\Gamma })=[-1. 1]$ 
for the singular potential $W(r)=W_{\Gamma }(r)=I_{[-1, 1]}(r)-r^{2}/2$. 
~\\
$\bullet $ $\beta (r)=\beta _{\Gamma }(r)=\partial I_{[-1, 1]}(r)+r^{3}$, 
$\pi (r)=\pi _{\Gamma }(r)=-r$ for all $r\in D(\beta )=D(\beta _{\Gamma })=[-1. 1]$ 
for the modified prototype double well potential $W(r)=W_{\Gamma }(r)=I_{[-1, 1]}(r)+(1/4)(r^{2}-1)^{2}-r^{2}/2$. \\
~\\
The third example is modified due to the assumption (A5). 
However, we can choose the original prototype double well potential (see Remark 4.1).

Our main theorem is given now. 

\paragraph{Theorem 2.1.}
{\it Under the assumptions {\rm (A1)}--{\rm (A4)}, for any given $m_{0}\in {\rm int} D(\beta _{\Gamma })$, 
there exists a weak solution of {\rm (P)} such that }
$$\frac{1}{|\Omega |+|\Gamma |}\left\{\int_{\Omega }u(0)dx+\int_{\Gamma }u_{\Gamma }(0)\right\}=m_{0}. $$

%%%%% Section 3. %%%%%
\section{Approximate problem and uniform estimates}
\setcounter{equation}{0}

In section, we consider the approximate problem 
and obtain the uniform estimates to show the existence of weak solutions of (P).

\subsection{Abstract formulation}
\indent

In order to prove the main theorem, we apply the abstract theory of evolution equation.
To do so, we define a proper lower semicontinuous convex functional $\varphi : \boldsymbol{H}_{0}\to [0, +\infty ]$ by
\begin{equation*}
	\label{varp}
	\varphi (\boldsymbol{z}):=\begin{cases}
	\displaystyle \frac{\kappa _{1}}{2}\int_{\Omega }|\nabla z|^{2}dx 
	+\frac{\kappa _{2}}{2}\int_{\Gamma }|\nabla _{\Gamma }z_{\Gamma }|^{2}d\Gamma \\
	~~~~~~~~~~~~\displaystyle +\int_{\Omega }\widehat{\beta }(z+m_{0})dx 
	+\int_{\Gamma }\widehat{\beta }_{\Gamma }(z_{\Gamma }+m_{0})d\Gamma \\
	~~~~~~~~~~~~~~~~~~~~ {\rm if} 
	\quad \boldsymbol{z}\in \boldsymbol{V}_{0} ~{\rm with~}\widehat{\beta }(z+m_{0})\in L^{1}(\Omega ), \quad \widehat{\beta }_{\Gamma }(z_{\Gamma }+m_{0})\in L^{1}(\Gamma ), \\
	+\infty  ~~~~~~~~~~~~~~~{\rm otherwise}. \\
\end{cases} 
\end{equation*}

Next, for each $\varepsilon \in (0, 1]$, we define a proper lower semicontinuous convex functional $\varphi _{\varepsilon }: \boldsymbol{H}_{0}\to [0, +\infty ]$ by
\begin{equation*}
	\varphi _{\varepsilon }(\boldsymbol{z}):=\begin{cases}
	\displaystyle \frac{\kappa _{1}}{2}\int_{\Omega }|\nabla z|^{2}dx 
	+\frac{\kappa _{2}}{2}\int_{\Gamma }|\nabla _{\Gamma }z_{\Gamma }|^{2}d\Gamma \\
	~~~~~~~~~~~~\displaystyle +\int_{\Omega }\widehat{\beta }_{\varepsilon }(z+m_{0})dx 
	+\int_{\Gamma }\widehat{\beta }_{\Gamma , \varepsilon }(z_{\Gamma }+m_{0})d\Gamma \\
	~~~~~~~~~~~~~~~~~~~~ {\rm if} 
	\quad \boldsymbol{z}\in \boldsymbol{V}_{0} ~{\rm with~}\widehat{\beta }_{\varepsilon }(z+m_{0})\in L^{1}(\Omega ), \quad \widehat{\beta }_{\Gamma , \varepsilon }(z_{\Gamma }+m_{0})\in L^{1}(\Gamma ), \\
	+\infty  ~~~~~~~~~~~~~~~{\rm otherwise}, \\
\end{cases} 
\end{equation*}
where $\widehat{\beta }_{\varepsilon }, \widehat{\beta }_{\Gamma , \varepsilon }$ defined as follows are {M}oreau--Yosida regularizations of $\widehat{\beta }, \widehat{\beta _{\Gamma }}$, respectively: 
$$\widehat{\beta }_{\varepsilon }(r):=\inf_{s\in \mathbb{R}}\left\{\frac{1}{2\varepsilon }|r-s|^{2}+\widehat{\beta }(s)\right\}=\frac{1}{2\varepsilon }|r-J_{\varepsilon }(r)|^{2}+\widehat{\beta }\bigl(J_{\varepsilon }(r)\bigr), $$
$$\widehat{\beta }_{\Gamma , \varepsilon }:=\inf_{s\in \mathbb{R}}\left\{\frac{1}{2\varepsilon \rho }|r-s|^{2}+\widehat{\beta }_{\Gamma }(s)\right\}=\frac{1}{2\varepsilon \rho }|r-J_{\Gamma , \varepsilon }(r)|^{2}+\widehat{\beta }_{\Gamma }\bigl(J_{\Gamma , \varepsilon }(r)\bigr), $$
where $\rho $ is a constant as in (\ref{ccbeta}) and $J_{\varepsilon }, J_{\Gamma , \varepsilon }: \mathbb{R}\to \mathbb{R}$ is resolvent operator given by 
$$J_{\varepsilon }(r):=(I+\varepsilon \beta )^{-1}(r), \quad J_{\Gamma , \varepsilon }(r):=(I+\varepsilon \rho \beta _{\Gamma })^{-1}(r)$$
for all $r\in \mathbb{R}$. 
Moreover, $\beta _{\varepsilon }, \beta _{\Gamma , \varepsilon }: \mathbb{R}\to \mathbb{R}$ defined as follows are Yosida approximatioon for maximal monotone operators $\beta , \beta _{\Gamma }$, respectively: 
$$\beta _{\varepsilon }(r):=\frac{1}{\varepsilon }\bigl(r-J_{\varepsilon }(r)\bigr), \quad \beta _{\Gamma , \varepsilon }(r):=\frac{1}{\varepsilon \rho }\bigl(r-J_{\Gamma , \varepsilon }(r)\bigr)$$
for all $r\in \mathbb{R}$, where $J_{\varepsilon }, J_{\Gamma , \varepsilon }: \mathbb{R}\to \mathbb{R}$ are resolvent operators. 
Then, it holds $\beta _{\varepsilon }(0)=\beta _{\Gamma , \varepsilon }(0)=0$. 
It is well known that $\beta _{\varepsilon }, \beta _{\Gamma , \varepsilon }$ are {L}ipschitz continuous 
with {L}ipschitz constant $1/\varepsilon , 1/(\varepsilon \rho )$, respectively. 
Here, we have following properties: 
$$0\leq \widehat{\beta }_{\varepsilon }(r)\leq \widehat{\beta }(r), \quad 0\leq \widehat{\beta }_{\Gamma , \varepsilon }(r)\leq \widehat{\beta }_{\Gamma }(r) \quad {\rm for~all~}r\in \mathbb{R}.$$
Hence, it holds $0\leq \varphi _{\varepsilon }(\boldsymbol{z})\leq \varphi (\boldsymbol{z})$ for all $\boldsymbol{z}\in \boldsymbol{H}_{0}$. 
These properties of {Y}osida approximation and {M}oreau--Yosida regularizations are as in \cite{Bar10, Bre73, Ken07}. 
Moreover, thanks to \cite[Lemma 4.4]{CC13}, it holds 
\begin{equation}\label{betaine}
|\beta _{\varepsilon }(r)|\leq \rho |\beta _{\Gamma , \varepsilon }(r)|+c_{0} \quad {\rm for~all}~r\in \mathbb{R} 
\end{equation}
with the same constants $\rho $ and $c_{0}$ as in (\ref{ccbeta}). 

Now, for each $\varepsilon \in (0, 1]$, we also define two proper lower semicontinuous convex functionals $\widetilde{\varphi }, \psi _{\varepsilon }: \boldsymbol{H}_{0}\to [0, +\infty ]$ by
\begin{equation*}
	\widetilde{\varphi }(\boldsymbol{z}):=\begin{cases}
	\displaystyle \frac{\kappa _{1}}{2}\int_{\Omega }|\nabla z|^{2}dx 
	+\frac{\kappa _{2}}{2}\int_{\Gamma }|\nabla _{\Gamma }z_{\Gamma }|^{2}d\Gamma & {\rm if} 
	\quad \boldsymbol{z}\in \boldsymbol{V}_{0}, \\
	+\infty  & {\rm otherwise} 
\end{cases} 
\end{equation*}
and 
\begin{equation*}
	\psi _{\varepsilon }(\boldsymbol{z}):=\begin{cases}
	\displaystyle \int_{\Omega }\widehat{\beta }_{\varepsilon }(z+m_{0})dx 
	+\int_{\Gamma }\widehat{\beta }_{\Gamma , \varepsilon }(z_{\Gamma }+m_{0})d\Gamma & {\rm if} 
	\quad \boldsymbol{z}\in \boldsymbol{V}_{0}, \\ %~{\rm with~}\widehat{\beta }_{\varepsilon }(z+m_{0})\in L^{1}(\Omega ), \quad \widehat{\beta }_{\Gamma , \varepsilon }(z_{\Gamma }+m_{0})\in L^{1}(\Gamma ), \\
	+\infty & {\rm otherwise}. 
\end{cases} 
\end{equation*}
Then, from \cite[Lemma C]{CF15}, the subdifferential $\boldsymbol{A}:=\partial _{\boldsymbol{H}_{0}}\widetilde{\varphi }$ on $\boldsymbol{H}_{0}$ is characterized by 
\begin{equation*}
\boldsymbol{A}(\boldsymbol{z})=(-\Delta z, \partial _{\boldsymbol{\nu }}z-\Delta _{\Gamma }z_{\Gamma }) 
\quad {\rm with}~\boldsymbol{z}=(z, z_{\Gamma })\in D(\boldsymbol{A})=\boldsymbol{W}\cap \boldsymbol{V}_{0}. 
\end{equation*}
Moreover, the representation of the subdifferential $\partial _{\boldsymbol{H}_{0}}\psi _{\varepsilon }$ is given by 
\begin{equation*}\label{subpsi}
\partial _{\boldsymbol{H}_{0}}\psi _{\varepsilon }(\boldsymbol{z})=\boldsymbol{P}\boldsymbol{\beta }_{\varepsilon }(\boldsymbol{z}+m_{0}\boldsymbol{1})
\quad {\rm for~all~}\boldsymbol{z}\in \boldsymbol{H}_{0}. 
\end{equation*}
This is proved by the same way as \cite[Lemma 3.2]{FM18}. 
Noting that it holds $D(\partial _{\boldsymbol{H}_{0}}\psi _{\varepsilon })=\boldsymbol{H}_{0}$ 
and $\boldsymbol{A}$ is a maximal monotone operator, 
it follows from the abstract monotonicity methods (see e.g., \cite[Sect.\ 2.1]{Bar10}) that 
$\boldsymbol{A}+\partial _{\boldsymbol{H}_{0}}\psi _{\varepsilon }$ is also a maximal monotone operator. 
Moreover, by a simple calculation, we deduce that $(\boldsymbol{A}+\partial _{\boldsymbol{H}_{0}}\psi _{\varepsilon })\subset \partial _{\boldsymbol{H}_{0}}\varphi _{\varepsilon }$. 
Hence, it follows that, for any $\boldsymbol{z}\in \boldsymbol{H}_{0}$, 
\begin{equation}\label{subvarp}
\partial _{\boldsymbol{H}_{0}}\varphi _{\varepsilon }(\boldsymbol{z})
=\bigl(\boldsymbol{A}+\partial _{\boldsymbol{H}_{0}}\psi _{\varepsilon }\bigr)(\boldsymbol{z}) 
\end{equation}
(see e.g., \cite{DK99}). 
%(\ref{subpsi}) and (\ref{subvarp}) are proved precisely in the {A}ppendix. 

\subsection{Approximate problem for (P)}
\indent

Now, we consider the following approximate problem, say (P)$_{\varepsilon }$: 
for each $\varepsilon \in (0, 1]$ find $\boldsymbol{v}_{\varepsilon }:=(v_{\varepsilon }, v_{\Gamma , \varepsilon })$ satisfying
\begin{eqnarray}
&&\varepsilon \boldsymbol{v}_{\varepsilon }'(t)+\boldsymbol{F}^{-1}\boldsymbol{v}_{\varepsilon }'(t)
+\partial _{\boldsymbol{H}_{0}}\varphi _{\varepsilon }\bigl(\boldsymbol{v}_{\varepsilon }(t)\bigr) \nonumber \\
&&~~~~~~~~~+\boldsymbol{P}\bigl(\widetilde{\boldsymbol{\pi }}\bigl(\boldsymbol{v}_{\varepsilon }(t)+m_{0}\boldsymbol{1}\bigr)\bigr)
=\boldsymbol{Pf}(t) \quad {\rm in~}\boldsymbol{H}_{0} \quad {\rm for~a.a.~}t\in (0, T) \label{apo1}\\
&&\boldsymbol{v}_{\varepsilon }(0)=\boldsymbol{v}_{\varepsilon }(T) \quad {\rm in~}\boldsymbol{H}_{0}. \label{apo2}
\end{eqnarray}
where, for all $\boldsymbol{z}\in \boldsymbol{H}$, 
%$\boldsymbol{\beta }_{\varepsilon }(\boldsymbol{z}):=(\beta _{\varepsilon }(z), \beta _{\Gamma , \varepsilon }(z_{\Gamma }))$, 
$\widetilde{\boldsymbol{\pi }}(\boldsymbol{z}):=(\widetilde{\pi }(z), \widetilde{\pi }_{\Gamma }(z_{\Gamma }))$ 
is cut-off function of $\pi , \pi _{\Gamma }$ given by 
\begin{eqnarray}\label{pi}
\widetilde{\pi }(r):=\begin{cases}
0 & {\rm if~}r\leq \sigma _{*}-1, \\
\pi (\sigma _{*})(r-\sigma _{*}+1) & {\rm if~}\sigma _{*}-1\leq r\leq \sigma _{*}, \\
\pi (r) & {\rm if~}\sigma _{*}\leq r\leq \sigma ^{*}, \\
-\pi (\sigma ^{*})(r-\sigma ^{*}-1) & {\rm if~}\sigma ^{*}\leq r\leq \sigma ^{*}+1, \\
0 & {\rm if~}r\geq \sigma ^{*}+1 
\end{cases}
\end{eqnarray}
and 
\begin{eqnarray}\label{pig}
\widetilde{\pi }_{\Gamma }(r):=\begin{cases}
0 & {\rm if~}r\leq \sigma _{\Gamma *}-1, \\
\pi _{\Gamma }(\sigma _{\Gamma *})(r-\sigma _{*}+1) & {\rm if~}\sigma _{\Gamma *}-1\leq r\leq \sigma _{\Gamma *}, \\
\pi _{\Gamma }(r) & {\rm if~}\sigma _{\Gamma *}\leq r\leq \sigma _{\Gamma }^{*}, \\
-\pi _{\Gamma }(\sigma _{\Gamma }^{*})(r-\sigma ^{*}-1) & {\rm if~}\sigma _{\Gamma }^{*}\leq r\leq \sigma _{\Gamma }^{*}+1, \\
0 & {\rm if~}r\geq \sigma _{\Gamma }^{*}+1 
\end{cases}
\end{eqnarray}
for all $r\in \mathbb{R}$, respectively. 
We establish the above cut-off function by referring to \cite{WZ16}.

From now, we show the next proposition. 

\paragraph{Proposition 3.1.}
{\it Under the assumptions {\rm (A1)}--{\rm (A5)}, for each $\varepsilon \in (0, 1]$, there exists a unique function 
$$\boldsymbol{v}_{\varepsilon }\in H^{1}(0, T; \boldsymbol{H}_{0})\cap L^{\infty }(0, T; \boldsymbol{V}_{0})\cap L^{2}(0, T; \boldsymbol{W})$$
such that $\boldsymbol{v}_{\varepsilon }$ satisfies {\rm (\ref{apo1})} and {\rm (\ref{apo2})}. } \\

In order to show the Proposition 3.1, 
we use the method in \cite{WZ16}, 
that is, we employ the viscosity approach. 
Conforming to the method, we consider the following problem: 
for each $\varepsilon \in (0, 1]$ and 
$\widetilde{\boldsymbol{f}}\in L^{2}(0, T; \boldsymbol{V}_{0})$, 
\begin{eqnarray}
&&(\boldsymbol{F}^{-1}+\varepsilon I)\boldsymbol{v}_{\varepsilon }'(t)+\partial \varphi _{\varepsilon }\bigl(\boldsymbol{v}_{\varepsilon }(t)\bigr)=\widetilde{\boldsymbol{f}}(t) \quad {\rm in~}\boldsymbol{H}_{0} \quad {\rm for~a.a.~}t\in (0, T), \label{ap2-1} \\
&&\boldsymbol{v}_{\varepsilon }(0)=\boldsymbol{v}_{\varepsilon }(T) \quad {\rm in~}\boldsymbol{H}_{0}. \label{ap2-2} 
\end{eqnarray}

Now, we can apply the abstract theory of doubly nonlinear evolution equation respect to time periodic problem \cite{AS11} for (\ref{ap2-1}), (\ref{ap2-2}) 
because the operator $\varepsilon \boldsymbol{I}+\boldsymbol{F}^{-1}$ 
and $\partial \varphi _{\varepsilon }$ are coercive in $\boldsymbol{H}_{0}$. 
%and $\varphi $ fulfills the stronger coercivity. 
This is an important assumption to use Theorem 2.2 in \cite{AS11}. 
It is an advantage of the viscosity approach. 
Hence, we obtain the next proposition.

\paragraph{Proposition 3.2.}
{\it For each $\varepsilon \in (0, 1]$ and $\widetilde{\boldsymbol{f}}\in L^{2}(0, T; \boldsymbol{V}_{0})$, 
there exists a unique function $\boldsymbol{v}_{\varepsilon }$ such that {\rm (\ref{ap2-1})} and {\rm (\ref{ap2-2})}. } \\

Hereafter, we apply the {S}chauder fixed point theorem to prove 
the existence of the solution of the problem (P)$_{\varepsilon }$. 
To this aim, we set 
$$\boldsymbol{Y}_{1}:=\bigl\{\bar{\boldsymbol{v}}_{\varepsilon }
\in H^{1}(0, T; \boldsymbol{H}_{0})\cap L^{\infty }(0, T; \boldsymbol{V}_{0}); 
\bar{\boldsymbol{v}}_{\varepsilon }(0)=\bar{\boldsymbol{v}}_{\varepsilon }(T)\bigr\}. $$
Firstly, for each $\bar{\boldsymbol{v}}_{\varepsilon }\in \boldsymbol{Y}_{1}$, 
we consider the following problem, say (P$_{\varepsilon }; \bar{\boldsymbol{v}}_{\varepsilon })$: 
\begin{equation}
\varepsilon \boldsymbol{v}_{\varepsilon }'(s)+\boldsymbol{F}^{-1}\boldsymbol{v}_{\varepsilon }'(s)
+\partial \varphi _{\varepsilon }\bigl(\boldsymbol{v}_{\varepsilon }(s)\bigr)
+\boldsymbol{P}\bigl(\widetilde{\boldsymbol{\pi }}\bigl(\bar{\boldsymbol{v}}_{\varepsilon }(s)+m_{0}\boldsymbol{1}\bigr)\bigr)
=\boldsymbol{Pf}(s) \quad {\rm in~}\boldsymbol{H}_{0} \label{Pa1}
\end{equation}
for a.a.\ $s\in (0, T)$, with
\begin{equation*}
\boldsymbol{v}_{\varepsilon }(0)=\boldsymbol{v}_{\varepsilon }(T) \quad {\rm in~}\boldsymbol{H}_{0}. \label{Pa2}
\end{equation*}

Next, we obtain the estimates of the solution of (P$_{\varepsilon }; \bar{\boldsymbol{v}}_{\varepsilon })$ 
to apply the {S}chauder fixed point theorem. 
Note that we can allow the dependent of $\varepsilon \in (0, 1]$ for estimates of Lemma 3.1 
because we use the {S}chauder fixed point theorem in the level of approximation. 

\paragraph{Lemma 3.1.}
{\it Let $\boldsymbol{v}_{\varepsilon }$ be the solution of problem {\rm (P}$_{\varepsilon }; \bar{\boldsymbol{v}}_{\varepsilon }{\rm )}$, 
it holds the following estimates. 
There exist positive constants 
$C_{1\varepsilon }, C_{2}, C_{3\varepsilon }$ such that 
\begin{gather}
\varepsilon \int_{0}^{T}\bigl|\boldsymbol{v}_{\varepsilon }'(s)\bigr|_{\boldsymbol{H}_{0}}^{2}ds
+\int_{0}^{T}\bigl|\boldsymbol{v}_{\varepsilon }'(s)\bigr|_{\boldsymbol{V}_{0}^{*}}^{2}ds
\leq C_{1\varepsilon }, \label{es3.1-1} \\
%\int_{0}^{T}\bigl|\boldsymbol{f}(s)\bigr|_{\boldsymbol{V}}^{2}ds+\frac{MT}{\varepsilon }, $$
\int_{0}^{T}\bigl|\boldsymbol{v}_{\varepsilon }(s)\bigr|_{\boldsymbol{V}_{0}}^{2}ds
+\int_{0}^{T}\int_{\Omega }\widehat{\beta }_{\varepsilon }\bigl(v_{\varepsilon }(s)+m_{0}\bigr)dxds
+\int_{0}^{T}\int_{\Gamma }\widehat{\beta }_{\Gamma , \varepsilon }\bigl(v_{\Gamma , \varepsilon }(s)+m_{0}\bigr)ds\leq C_{2} \label{es3.1-2}
\end{gather}
and 
\begin{gather}
\frac{1}{2}\bigl|\boldsymbol{v}_{\varepsilon }(t)\bigr|_{\boldsymbol{V}_{0}}^{2}
+\int_{\Omega }\widehat{\beta }_{\varepsilon }\bigl(v_{\varepsilon }(t)+m_{0}\bigr)dx
+\int_{\Gamma }\widehat{\beta }_{\Gamma , \varepsilon }\bigl(v_{\Gamma , \varepsilon }(t)+m_{0}\bigr)d\Gamma 
\leq C_{3\varepsilon } \label{es3.1-3}
\end{gather}
for a.a.\ $t\in (0, T)$. }

\paragraph{Proof}
At first, for each $\bar{\boldsymbol{v}}_{\varepsilon }\in Y_{1}$, note that there exists a positive constant $M$, depending only on 
$\sigma _{*}, \sigma _{\Gamma *}, \sigma ^{*}$ and $\sigma _{\Gamma }^{*}$, such that 
\begin{equation}\label{pibdd}
\bigl|\widetilde{\boldsymbol{\pi }}\bigl(\bar{\boldsymbol{v}}_{\varepsilon }(t)
+m_{0}\boldsymbol{1}\bigr)\bigr|_{\boldsymbol{H}_{0}}^{2}
\leq M \quad {\rm for~all~}t\in [0, T]. 
\end{equation}
Now, testing (\ref{Pa1}) at time $s\in (0, T)$ by $\boldsymbol{v}_{\varepsilon }'(s)$ 
and using the {Y}oung inequality, we infer that 
\begin{eqnarray*}
&&\varepsilon \bigl|\boldsymbol{v}_{\varepsilon }'(s)\bigr|_{\boldsymbol{H}_{0}}^{2}
+\bigl|\boldsymbol{v}_{\varepsilon }'(s)\bigr|_{\boldsymbol{V}_{0}^{*}}^{2}
+\frac{d}{ds}\varphi _{\varepsilon }\bigl(\boldsymbol{v}_{\varepsilon }(s)\bigr) \\
&&\quad =\bigl(\boldsymbol{Pf}(s)-
\boldsymbol{P}\bigl(\widetilde{\boldsymbol{\pi }}\bigl(\bar{\boldsymbol{v}}_{\varepsilon }(s)+m_{0}\boldsymbol{1}\bigr)\bigr), 
\boldsymbol{v}_{\varepsilon }'(s)\bigr)_{\boldsymbol{H}_{0}} \\
&&\quad \leq \frac{1}{2}\bigl|\boldsymbol{f}(s)\bigr|_{\boldsymbol{V}}^{2}
+\frac{1}{2}\bigl|\boldsymbol{v}_{\varepsilon }'(s)\bigr|_{\boldsymbol{V}_{0}^{*}}^{2}
+\frac{M}{2\varepsilon }
+\frac{\varepsilon }{2}\bigl|\boldsymbol{v}_{\varepsilon }'(s)\bigr|_{\boldsymbol{H}_{0}}^{2} 
\end{eqnarray*}
for a.a.\ $s\in (0, T)$. 
Therefore, we have that 
\begin{equation}\label{ine2}
\varepsilon \bigl|\boldsymbol{v}_{\varepsilon }'(s)\bigr|_{\boldsymbol{H}_{0}}^{2}
+\bigl|\boldsymbol{v}_{\varepsilon }'(s)\bigr|_{\boldsymbol{V}_{0}^{*}}^{2}
+2\frac{d}{ds}\varphi _{\varepsilon }\bigl(\boldsymbol{v}_{\varepsilon }(s)\bigr)
\leq \bigl|\boldsymbol{f}(s)\bigr|_{\boldsymbol{V}}^{2}+\frac{M}{\varepsilon }. 
\end{equation}
Then, integrating it over $(0, T)$ with respect to $s$ and 
using the periodic property, we see that 
$$\varepsilon \int_{0}^{T}\bigl|\boldsymbol{v}_{\varepsilon }'(s)\bigr|_{\boldsymbol{H}_{0}}^{2}ds
+\int_{0}^{T}\bigl|\boldsymbol{v}_{\varepsilon }'(s)\bigr|_{\boldsymbol{V}_{0}^{*}}^{2}ds
\leq \int_{0}^{T}\bigl|\boldsymbol{f}(s)\bigr|_{\boldsymbol{V}}^{2}ds+\frac{MT}{\varepsilon }, $$
which implies that it follows the first estimate (\ref{es3.1-1}). 

On the other hand, 
testing (\ref{Pa1}) at time $s\in (0, T)$ by $\boldsymbol{v}_{\varepsilon }(s)$
and from (\ref{PW0}), we deduce that 
\begin{eqnarray*}
&&\frac{1}{2}\frac{d}{ds}\bigl|\boldsymbol{v}_{\varepsilon }(s)\bigr|_{\boldsymbol{V}_{0}^{*}}^{2}
+\frac{\varepsilon }{2}\frac{d}{ds}\bigl|\boldsymbol{v}_{\varepsilon }(s)\bigr|_{\boldsymbol{H}_{0}}^{2}
+\varphi _{\varepsilon }\bigl(\boldsymbol{v}_{\varepsilon }(s)\bigr) \\
&&\quad \leq \bigl(\boldsymbol{Pf}(s)-
\boldsymbol{P}\bigl(\widetilde{\boldsymbol{\pi }}\bigl(\bar{\boldsymbol{v}}_{\varepsilon }(s)+m_{0}\boldsymbol{1}\bigr)\bigr), 
\boldsymbol{v}_{\varepsilon }(s)\bigr)_{\boldsymbol{H}_{0}}+\varphi _{\varepsilon }(0) \\
&&\quad \leq 2c_{{\rm P}}\bigl|\boldsymbol{f}(s)\bigr|_{\boldsymbol{H}_{0}}^{2}+
\frac{1}{4c_{{\rm P}}}\bigl|\boldsymbol{v}_{\varepsilon }(s)\bigr|_{\boldsymbol{H}_{0}}^{2}+2c_{{\rm P}}M+\varphi (0) \\
&&\quad \leq 2c_{{\rm P}}\bigl|\boldsymbol{f}(s)\bigr|_{\boldsymbol{H}_{0}}^{2}+
\frac{1}{4}\bigl|\boldsymbol{v}_{\varepsilon }(s)\bigr|_{\boldsymbol{V}_{0}}^{2}+2c_{{\rm P}}M+\varphi (0) \\
&&\quad \leq 2c_{{\rm P}}\bigl|\boldsymbol{f}(s)\bigr|_{\boldsymbol{H}_{0}}^{2}+
\frac{1}{2}\varphi _{\varepsilon }\bigl(\boldsymbol{v}_{\varepsilon }(s)\bigr)+2c_{{\rm P}}M+\varphi (0)
\end{eqnarray*}
for a.a.\ $s\in (0, T)$, thanks to the definition of the subdifferential.  
From the definition of $\varphi _{\varepsilon }$, it follows that 
\begin{eqnarray*}
&&\frac{1}{2}\frac{d}{ds}\bigl|\boldsymbol{v}_{\varepsilon }(s)\bigr|_{\boldsymbol{V}_{0}^{*}}^{2}
+\frac{\varepsilon }{2}\frac{d}{ds}\bigl|\boldsymbol{v}_{\varepsilon }(s)\bigr|_{\boldsymbol{H}_{0}}^{2}
+\frac{1}{4}\bigl|\boldsymbol{v}_{\varepsilon }(s)\bigr|_{\boldsymbol{V}_{0}}^{2} \\
&&\quad \quad \quad +\frac{1}{2}\int_{\Omega }\widehat{\beta }_{\varepsilon }\bigl(v_{\varepsilon }(s)+m_{0}\bigr)dx
+\frac{1}{2}\int_{\Gamma }\widehat{\beta }_{\Gamma , \varepsilon }\bigl(v_{\Gamma , \varepsilon }(s)+m_{0}\bigr)d\Gamma \\
&&\quad \leq 2c_{{\rm P}}\bigl|\boldsymbol{f}(s)\bigr|_{\boldsymbol{H}_{0}}^{2}+2c_{{\rm P}}M+
+\int_{\Omega }\widehat{\beta }_{\varepsilon }\bigl(m_{0}\bigr)dx
+\int_{\Gamma }\widehat{\beta }_{\Gamma , \varepsilon }\bigl(m_{0}\bigr)d\Gamma 
\end{eqnarray*}
for a.a.\ $s\in (0, T)$. 
Integrating it over $(0, T)$ and using the periodic property, we see that 
\begin{eqnarray*}
&&\frac{1}{2}\int_{0}^{T}\bigl|\boldsymbol{v}_{\varepsilon }(s)\bigr|_{\boldsymbol{V}_{0}}^{2}ds
+\int_{0}^{T}\int_{\Omega }\widehat{\beta }_{\varepsilon }\bigl(v_{\varepsilon }(s)+m_{0}\bigr)dxds
+\int_{0}^{T}\int_{\Gamma }\widehat{\beta }_{\Gamma , \varepsilon }\bigl(v_{\Gamma , \varepsilon }(s)+m_{0}\bigr)d\Gamma ds\\
&&\quad \leq 4c_{{\rm P}}\bigl|\boldsymbol{f}\bigr|_{L^{2}(0, T; \boldsymbol{H}_{0})}^{2}
+4c_{{\rm P}}TM
+T|\Omega |\bigl|\widehat{\beta }\bigl(m_{0}\bigr)\bigr|
+T|\Gamma |\bigl|\widehat{\beta }_{\Gamma }\bigl(m_{0}\bigr)\bigr|. 
\end{eqnarray*}
Hence, there exist a positive constant $C_{2}$ such that the second estimate (\ref{es3.1-2}) holds. 

Next, for each $s, t\in [0, T]$ such that $s\leq t$, 
we integrate (\ref{ine2}) over $[s, t]$ with respect to $s$. 
Then, by neglecting the first two positive terms, we have 
$$\varphi _{\varepsilon }\bigl(\boldsymbol{v}_{\varepsilon }(t)\bigr)
\leq \varphi _{\varepsilon }\bigl(\boldsymbol{v}_{\varepsilon }(s)\bigr)
+\frac{1}{2}\int_{0}^{T}\bigl|\boldsymbol{f}(s)\bigr|_{\boldsymbol{V}}^{2}ds
+\frac{MT}{2\varepsilon }$$
for a.a.\ $s, t\in [0, T]$, 
namely, 
\begin{eqnarray}
&&\frac{1}{2}\bigl|\boldsymbol{v}_{\varepsilon }(t)\bigr|_{\boldsymbol{V}_{0}}^{2}
+\int_{\Omega }\widehat{\beta }_{\varepsilon }\bigl(v_{\varepsilon }(t)+m_{0}\bigr)dx
+\int_{\Gamma }\widehat{\beta }_{\Gamma , \varepsilon }\bigl(v_{\Gamma , \varepsilon }(t)+m_{0}\bigr)d\Gamma \nonumber \\
&&\quad \leq \frac{1}{2}\bigl|\boldsymbol{v}_{\varepsilon }(s)\bigr|_{\boldsymbol{V}_{0}}^{2}
+\int_{\Omega }\widehat{\beta }_{\varepsilon }\bigl(v_{\varepsilon }(s)+m_{0}\bigr)dx
+\int_{\Gamma }\widehat{\beta }_{\Gamma , \varepsilon }\bigl(v_{\Gamma , \varepsilon }(s)+m_{0}\bigr)d\Gamma \nonumber \\
&&\quad \quad \quad +\frac{1}{2}\int_{0}^{T}\bigl|\boldsymbol{f}(s)\bigr|_{\boldsymbol{V}}^{2}ds
+\frac{MT}{2\varepsilon }  \label{st}
\end{eqnarray}
for a.a.\ $s, t\in [0, T]$. 
Now, integrating it over $(0, t)$ with respect to $s$, we deduce that 
\begin{eqnarray}
&&\frac{t}{2}\bigl|\boldsymbol{v}_{\varepsilon }(t)\bigr|_{\boldsymbol{V}_{0}}^{2}
+t\int_{\Omega }\widehat{\beta }_{\varepsilon }\bigl(v_{\varepsilon }(t)+m_{0}\bigr)dx
+t\int_{\Gamma }\widehat{\beta }_{\Gamma , \varepsilon }\bigl(v_{\Gamma , \varepsilon }(t)+m_{0}\bigr)d\Gamma \nonumber \\
&&\quad \leq \frac{1}{2}\int_{0}^{T}\bigl|\boldsymbol{v}_{\varepsilon }(s)\bigr|_{\boldsymbol{V}_{0}}^{2}ds
+\int_{0}^{T}\int_{\Omega }\widehat{\beta }_{\varepsilon }\bigl(v_{\varepsilon }(s)+m_{0}\bigr)dxds
+\int_{0}^{T}\int_{\Gamma }\widehat{\beta }_{\Gamma , \varepsilon }\bigl(v_{\Gamma , \varepsilon }(s)+m_{0}\bigr)d\Gamma ds \nonumber \\
&&\quad \quad \quad +\frac{T}{2}\int_{0}^{T}\bigl|\boldsymbol{f}(s)\bigr|_{\boldsymbol{V}}^{2}ds
+\frac{MT^{2}}{2\varepsilon } \label{allt}
\end{eqnarray}
for a.a.\ $t\in [0, T]$. 
In particular, putting $t:=T$ and dividing (\ref{allt}) by $T$, it follows that 
\begin{eqnarray}
&&\frac{1}{2}\bigl|\boldsymbol{v}_{\varepsilon }(T)\bigr|_{\boldsymbol{V}_{0}}^{2}
+\int_{\Omega }\widehat{\beta }_{\varepsilon }\bigl(v_{\varepsilon }(T)+m_{0}\bigr)dx
+\int_{\Gamma }\widehat{\beta }_{\Gamma , \varepsilon }\bigl(v_{\Gamma , \varepsilon }(T)+m_{0}\bigr)d\Gamma \nonumber \\
&&\quad \leq \frac{1}{2T}\int_{0}^{T}\bigl|\boldsymbol{v}_{\varepsilon }(s)\bigr|_{\boldsymbol{V}_{0}}^{2}ds
+\frac{1}{T}\int_{0}^{T}\int_{\Omega }\widehat{\beta }_{\varepsilon }\bigl(v_{\varepsilon }(s)+m_{0}\bigr)dxds \nonumber \\
&&\quad \quad \quad +\frac{1}{T}\int_{0}^{T}\int_{\Gamma }\widehat{\beta }_{\Gamma , \varepsilon }
\bigl(v_{\Gamma , \varepsilon }(s)+m_{0}\bigr)d\Gamma ds
+\frac{1}{2}\int_{0}^{T}\bigl|\boldsymbol{f}(s)\bigr|_{\boldsymbol{V}}^{2}ds
+\frac{MT}{2\varepsilon }. \label{ineT}
\end{eqnarray}
Hence, combining the second estimste (\ref{es3.1-2}) and (\ref{ineT}), we see that 
\begin{eqnarray*}
&&\frac{1}{2}\bigl|\boldsymbol{v}_{\varepsilon }(T)\bigr|_{\boldsymbol{V}_{0}}^{2}
+\int_{\Omega }\widehat{\beta }_{\varepsilon }\bigl(v_{\varepsilon }(T)+m_{0}\bigr)dx
+\int_{\Gamma }\widehat{\beta }_{\Gamma , \varepsilon }\bigl(v_{\Gamma , \varepsilon }(T)+m_{0}\bigr)d\Gamma \nonumber \\
&&\quad \leq \frac{C_{2}}{T}
+\frac{1}{2}\int_{0}^{T}\bigl|\boldsymbol{f}(s)\bigr|_{\boldsymbol{V}}^{2}ds
+\frac{MT}{2\varepsilon }. 
\end{eqnarray*}
Moreover, from the periodic property, we infer that 
\begin{eqnarray}
&&\frac{1}{2}\bigl|\boldsymbol{v}_{\varepsilon }(0)\bigr|_{\boldsymbol{V}_{0}}^{2}
+\int_{\Omega }\widehat{\beta }_{\varepsilon }\bigl(v_{\varepsilon }(0)+m_{0}\bigr)dx
+\int_{\Gamma }\widehat{\beta }_{\Gamma , \varepsilon }\bigl(v_{\Gamma , \varepsilon }(0)+m_{0}\bigr)d\Gamma \nonumber \\
&&\quad \leq \frac{C_{2}}{T}
+\frac{1}{2}\int_{0}^{T}\bigl|\boldsymbol{f}(s)\bigr|_{\boldsymbol{V}}^{2}ds
+\frac{MT}{2\varepsilon }. \label{es0}
\end{eqnarray}
Now, let $s$ be $0$ in (\ref{st}). 
Then, owing to (\ref{es0}), we deduce that 
\begin{eqnarray*}
&&\frac{1}{2}\bigl|\boldsymbol{v}_{\varepsilon }(t)\bigr|_{\boldsymbol{V}_{0}}^{2}
+\int_{\Omega }\widehat{\beta }_{\varepsilon }\bigl(v_{\varepsilon }(t)+m_{0}\bigr)dx
+\int_{\Gamma }\widehat{\beta }_{\Gamma , \varepsilon }\bigl(v_{\Gamma , \varepsilon }(t)+m_{0}\bigr)d\Gamma \\
%&&\quad \leq \bigl|\boldsymbol{v}_{\varepsilon }(0)\bigr|_{\boldsymbol{V}_{0}}^{2}
%+2\int_{\Omega }\widehat{\beta }\bigl(v_{\varepsilon }(0)+m_{0}\bigr)dx
%+2\int_{\Gamma }\widehat{\beta }_{\Gamma , \varepsilon }\bigl(v_{\Gamma , \varepsilon }(0)+m_{0}\bigr)d\Gamma \\
%&&\quad \quad \quad +\bigl|\boldsymbol{f}\bigr|_{L^{2}(0, T; \boldsymbol{V})}^{2}
%+\frac{MT}{\varepsilon } \\
&&\quad \leq \frac{C_{2}}{T}
+|\boldsymbol{f}|_{L^{2}(0, T; \boldsymbol{V})}^{2}
+\frac{MT}{\varepsilon }
\end{eqnarray*}
for a.a.\ $t\in [0, T]$. 
Therefore, there exisits a positive constant $C_{3\varepsilon }$ such that 
the final estimate (\ref{es3.1-3}) holds. 
\hfill $\Box $ \\

In terms of (\ref{es3.1-1}), the key point to prove the estimate is exploiting (\ref{pibdd}). 
(\ref{pibdd}) is arised from the form of cut-off functions (\ref{pi}), (\ref{pig}). 
The form of cut-off functions depend on the assumption (A5) essentially. 
However, considered the same estimate in \cite[Lemma 4.1]{CF15}, 
They do not impose the assumption. 
They use the {G}ronwall inequality to obtain the estimate 
because the initial value is given data. 
On the other hand, we can not obtain it even though we use the {G}ronwall inequality, 
becasuse the initial value is not given. 
For this reason, it is nacessary to impose (A5). 
This is a difficult point to solve the time periodic problem.

Now, we show the existence of solutions of approximate problem (P)$_{\varepsilon }$.

\paragraph{Proof of Proposition 3.1}
To this aim, 
we apply the {S}chauder fixed point theorem. 
To do so, we set 
$$\boldsymbol{Y}_{2}:=
\left\{\bar{\boldsymbol{v}}_{\varepsilon }\in \boldsymbol{Y}_{1}; 
\sup _{t\in [0, T]}\bigl|\bar{\boldsymbol{v}}_{\varepsilon }\bigr|_{\boldsymbol{V}_{0}}^{2}
+\varepsilon \bigl|\bar{\boldsymbol{v}}_{\varepsilon }\bigr|_{H^{1}(0, T; \boldsymbol{H}_{0})}^{2}
%+\bigl|\bar{\boldsymbol{v}}_{\varepsilon }\bigr|_{H^{1}(0, T; \boldsymbol{V}_{0}^{*})}^{2}
\leq M_{\varepsilon }\right\}, $$
where $M_{\varepsilon }$ is a positive constant and be determined by Lemma 3.1. 
Then, the set $\boldsymbol{Y}_{2}$ is non-empty compact convex on $C(0, T; \boldsymbol{H}_{0})$. 
Now, from Proposition 3.2, 
for each $\bar{\boldsymbol{v}}_{\varepsilon }\in \boldsymbol{Y}_{2}$, 
there exists a unique solution $\boldsymbol{v}_{\varepsilon }$ of (P$_{\varepsilon }; \bar{\boldsymbol{v}}_{\varepsilon })$. 
Moreover, from Lemma 3.1, it holds $\boldsymbol{v}_{\varepsilon }\in \boldsymbol{Y}_{2}$. 
Here, we define the mapping $\boldsymbol{S}: \boldsymbol{Y}_{2}\to \boldsymbol{Y}_{2}$ such that, 
for each $\bar{\boldsymbol{v}}_{\varepsilon }\in \boldsymbol{Y}_{2}$, 
corresponding $\bar{\boldsymbol{v}}_{\varepsilon }$ to 
the solution $\boldsymbol{v}_{\varepsilon }$ of (P$_{\varepsilon }; \bar{\boldsymbol{v}}_{\varepsilon })$. 
Then, the mapping $\boldsymbol{S}$ is continuous 
on $\boldsymbol{Y}_{2}$ with respect to topology of $C(0, T; \boldsymbol{H}_{0})$. 
Indeed, let $\{\bar{\boldsymbol{v}}_{\varepsilon , n}\}_{n\in \mathbb{N}}\subset \boldsymbol{Y}_{2}$ 
be $\bar{\boldsymbol{v}}_{\varepsilon , n}\to \bar{\boldsymbol{v}}_{\varepsilon }$ 
in $C(0, T; \boldsymbol{H}_{0})$ 
and $\{\boldsymbol{v}_{\varepsilon , n}\}_{n\in \mathbb{N}}$ 
be the sequence of the solution of (P$_{\varepsilon }; \bar{\boldsymbol{v}}_{\varepsilon , n})$. 
From Lemma 3.1, there exist a subsequence $\{n_{k}\}_{k\in \mathbb{N}}$, with $n_{k}\to \infty $ as $k\to \infty $, 
and $\boldsymbol{v}_{\varepsilon }\in H^{1}(0, T; \boldsymbol{H}_{0})\cap L^{\infty }(0, T; \boldsymbol{V}_{0})$ such that 
\begin{equation}\label{vep}
\boldsymbol{v}_{\varepsilon , n_{k}}\to \boldsymbol{v}_{\varepsilon } 
\quad {\rm weakly~star~in~}H^{1}(0, T; \boldsymbol{H}_{0})\cap L^{\infty }(0, T; \boldsymbol{V}_{0}). 
\end{equation}
Hence, from (\ref{vep}) and the {A}scoli--{A}lzela theorem (see e.g., \cite{Sim87}), 
there exists a subsequence (not relabeled) such that 
\begin{equation}\label{strong}
\boldsymbol{v}_{\varepsilon , n_{k}}\to \boldsymbol{v}_{\varepsilon } \quad {\rm in~}C([0, T]; \boldsymbol{H}_{0})
\end{equation}
as $k\to \infty $. 
Also, we have 
\begin{equation}
\boldsymbol{v}_{\varepsilon , n_{k}}'\to \boldsymbol{v}_{\varepsilon }' \quad {\rm weakly~in~}L^{2}(0, T; \boldsymbol{H}_{0}) \label{dds}
\end{equation}
as $k\to \infty $. 
Because we have $\boldsymbol{v}_{\varepsilon , n_{k}}(0)=\boldsymbol{v}_{\varepsilon , n_{k}}(T)$, 
it holds $\boldsymbol{v}_{\varepsilon }(0)=\boldsymbol{v}_{\varepsilon }(T)$ in $\boldsymbol{H}_{0}$. 
Hereafter, we show that this $\boldsymbol{v}_{\varepsilon }$ is the solution of (P$_{\varepsilon }; \bar{\boldsymbol{v}}_{\varepsilon })$. 
Since $\boldsymbol{v}_{\varepsilon , n_{k}}$ is the solution 
of (P$_{\varepsilon }; \bar{\boldsymbol{v}}_{\varepsilon , n_{k}})$, 
we see that 
\begin{eqnarray}\label{defsub}
&&\int_{0}^{T}\bigl(\boldsymbol{Pf}(s)-\boldsymbol{P}\bigl(\widetilde{\boldsymbol{\pi }}\bigl(\bar{\boldsymbol{v}}_{\varepsilon , n_{k}}(s)
+m_{0}\bigr)\bigr)-\varepsilon \boldsymbol{v}_{\varepsilon , n_{k}}'(s)
-\boldsymbol{F}^{-1}\boldsymbol{v}_{\varepsilon , n_{k}}'(s), 
\boldsymbol{\eta }(s)-\boldsymbol{v}_{\varepsilon , n_{k}}(s)\bigr)_{\boldsymbol{H}_{0}}ds \nonumber \\
&&\quad \leq \int_{0}^{T}\varphi _{\varepsilon }\bigl(\boldsymbol{\eta }(s)\bigr)ds
-\int_{0}^{T}\varphi _{\varepsilon }\bigl(\boldsymbol{v}_{\varepsilon , n_{k}}(s)\bigr)ds
\end{eqnarray}
for all $\boldsymbol{\eta }\in L^{2}(0, T; \boldsymbol{H}_{0})$, 
thanks to the definition of subdifferential $\partial \varphi _{\varepsilon }$. 
%Now, we have 
%\begin{equation}
%\int_{0}^{T}\frac{d}{ds}\bigl(\boldsymbol{F}^{-1}\bigl(\boldsymbol{v}_{\varepsilon , n_{k}}(s), 
%\boldsymbol{v}_{\varepsilon , n_{k}}(s)\bigr)\bigr)
%\end{equation}
%\begin{equation}
%\end{equation}
Moreover, it follows from 
$\bar{\boldsymbol{v}}_{\varepsilon , n_{k}}\to \bar{\boldsymbol{v}}_{\varepsilon }$ in $C(0, T; \boldsymbol{H}_{0})$ that 
\begin{equation}\label{ppi}
\boldsymbol{P}\bigl(\widetilde{\boldsymbol{\pi }}(\bar{\boldsymbol{v}}_{\varepsilon , n_{k}}+m_{0})\bigr)
\to \boldsymbol{P}\bigl(\widetilde{\boldsymbol{\pi }}(\bar{\boldsymbol{v}}_{\varepsilon }+m_{0})\bigr)
\quad {\rm in~}C([0, T]; \boldsymbol{H}_{0}). 
\end{equation}
Thus, on account of (\ref{vep})--(\ref{ppi}), 
by taking the upper limit as $k\to \infty $ in (\ref{defsub}) and using 
$$\liminf _{k\to \infty }\int_{0}^{T}\varphi _{\varepsilon }\bigl(\boldsymbol{v}_{\varepsilon , n_{k}}(s)\bigr)ds
\geq \int_{0}^{T}\varphi _{\varepsilon }\bigl(\boldsymbol{v}_{\varepsilon }(s)\bigr)ds, $$
we infer that 
\begin{eqnarray*}
&&\int_{0}^{T}\bigl(\boldsymbol{Pf}(s)
-\boldsymbol{P}\bigl(\widetilde{\boldsymbol{\pi }}\bigl(\bar{\boldsymbol{v}}_{\varepsilon }(s)+m_{0}\bigr)\bigr)
-\varepsilon \boldsymbol{v}_{\varepsilon }'(s)
-\boldsymbol{F}^{-1}\boldsymbol{v}_{\varepsilon }'(s), 
\boldsymbol{\eta }(s)
-\boldsymbol{v}_{\varepsilon }(s)\bigr)_{\boldsymbol{H}_{0}}ds \\
&&\quad \leq \int_{0}^{T}\varphi _{\varepsilon }\bigl(\boldsymbol{\eta }(s)\bigr)ds
-\int_{0}^{T}\varphi _{\varepsilon }\bigl(\boldsymbol{v}_{\varepsilon }(s)\bigr)ds
\end{eqnarray*}
for all $\boldsymbol{\eta }\in L^{2}(0, T; \boldsymbol{H}_{0})$. 
Hence, we see that the function $\boldsymbol{v}_{\varepsilon }$ 
is the solution of (P$_{\varepsilon }; \bar{\boldsymbol{v}}_{\varepsilon }$). 
As a result, 
it follows from the uniqueness of the solution of (P$_{\varepsilon }; \bar{\boldsymbol{v}}_{\varepsilon }$) that 
$$\boldsymbol{S}(\bar{\boldsymbol{v}}_{\varepsilon , n_{k}})
=\boldsymbol{v}_{\varepsilon , n_{k}}
\to \boldsymbol{v}_{\varepsilon }
=\boldsymbol{S}(\bar{\boldsymbol{v}}_{\varepsilon })
\quad {\rm in~}C([0, T]; \boldsymbol{H}_{0})$$
as $k\to \infty $. 
Therefore, 
the mapping $\boldsymbol{S}$ is continuous with respect to $C([0, T]; \boldsymbol{H}_{0})$. 
Thus, from the {S}chauder fixed point theorem, 
there exists a fixed point on $\boldsymbol{Y}_{2}$, 
namely, the problem (P)$_{\varepsilon }$ admits a solution $\boldsymbol{v}_{\varepsilon }$. 
Finally, from the fact that 
$\partial \varphi _{\varepsilon }(\boldsymbol{v}_{\varepsilon })\in L^{2}(0, T; \boldsymbol{H}_{0})$, 
which implies $\boldsymbol{v}_{\varepsilon }\in L^{2}(0, T; \boldsymbol{W})$. 
\hfill $\Box $ \\

%%%%mu'Ì'è‹`%%%%%%%%%%%
Now, 
%based on \cite{CF15}, 
we consider the chemical potential $\boldsymbol{\mu }:=(\mu , \mu _{\Gamma })$
by approximating. 
For each $\varepsilon \in (0, 1]$, we set the approximate sequence 
\begin{equation}\label{apromusub}
\boldsymbol{\mu }_{\varepsilon }(s):=
\varepsilon \boldsymbol{v}_{\varepsilon }'(s)
%+\boldsymbol{A}\boldsymbol{v}_{\varepsilon }(t)
%+\boldsymbol{\beta }_{\varepsilon }\bigl(\boldsymbol{v}_{\varepsilon }(t)+m_{0}\boldsymbol{1}\bigr)
+\partial \varphi _{\varepsilon }\bigl(\boldsymbol{v}_{\varepsilon }(s)\bigr)
+\widetilde{\boldsymbol{\pi }}\bigl(\boldsymbol{v}_{\varepsilon }(s)+m_{0}\boldsymbol{1}\bigr)
-\boldsymbol{f}(s)
\end{equation}
for a.a.\ $s\in (0, T)$. 
From (\ref{subvarp}), we can rewrite (\ref{apromusub}) as 
\begin{equation}\label{apromu}
\boldsymbol{\mu }_{\varepsilon }(s)=
\varepsilon \boldsymbol{v}_{\varepsilon }'(s)
+\boldsymbol{A}\boldsymbol{v}_{\varepsilon }(s)
+\boldsymbol{P}\boldsymbol{\beta }_{\varepsilon }\bigl(\boldsymbol{v}_{\varepsilon }(s)+m_{0}\boldsymbol{1}\bigr)
%+\partial \varphi _{\varepsilon }\bigl(\boldsymbol{v}_{\varepsilon }(t)\bigr)
+\widetilde{\boldsymbol{\pi }}\bigl(\boldsymbol{v}_{\varepsilon }(s)+m_{0}\boldsymbol{1}\bigr)
-\boldsymbol{f}(s)
\end{equation}
for a.a.\ $s\in (0, T)$. 
%This is the same setting for $\boldsymbol{\mu }_{\varepsilon }$ in \cite{CF15}. 
%Therefore, note that the result in terms of $\boldsymbol{\mu }_{\varepsilon }$ is the same \cite{CF15}. 
%Hereafter, since it is written particularly in \cite{CF15}, we only sketch the result about $\boldsymbol{\mu }_{\varepsilon }$. 
Then, %from the representation of the subdifferential operator $\partial \varphi _{\varepsilon }$, 
we rewrite (\ref{apo1}) as 
$$\boldsymbol{F}^{-1}\boldsymbol{v}_{\varepsilon }'(s)
+\boldsymbol{\mu }_{\varepsilon }(s)
-\omega _{\varepsilon }(s)\boldsymbol{1}
=\boldsymbol{0} \quad {\rm in~}\boldsymbol{V}$$
for a.a.\ $s\in (0, T)$, 
where 
$$\omega _{\varepsilon }(s):=
m\bigl(%\boldsymbol{\beta }_{\varepsilon }\bigl(\boldsymbol{v}_{\varepsilon }(t)+m_{0}\boldsymbol{1}\bigr)+
\widetilde{\boldsymbol{\pi }}\bigl(\boldsymbol{v}_{\varepsilon }(s)+m_{0}\boldsymbol{1}\bigr)
-\boldsymbol{f}(s)\bigr)$$
for a.a.\ $s\in (0, T)$. 
Therefore, we have 
$\boldsymbol{P}\boldsymbol{\mu }_{\varepsilon }=\boldsymbol{\mu }_{\varepsilon }-\omega _{\varepsilon }\boldsymbol{1}\in L^{2}(0, T; \boldsymbol{V}_{0})$ and 
$\omega _{\varepsilon }\in L^{2}(0, T)$. 
Then, it holds $\boldsymbol{\mu }_{\varepsilon }\in L^{2}(0, T; \boldsymbol{V})$ and 
\begin{equation}
\boldsymbol{v}_{\varepsilon }'(s)
+\boldsymbol{F}\boldsymbol{P}\boldsymbol{\mu }_{\varepsilon }(s)
=\boldsymbol{0} \quad {\rm in~}\boldsymbol{V}_{0}^{*} \label{vmu}
\end{equation}
for a.a.\ $s\in (0, T)$.

\subsection{Uniform estimates}
\indent

In this subsection, we obtain uniform estimates independent of $\varepsilon \in (0, 1]$. 
We refer to \cite{WZ16} to obtain uniform estimates.

%%%first estimate%%%%%%%%%%%%%%%%%%%%%%%%%%%%%%%%%%%%%%%%%%
\paragraph{Lemma 3.2.}
{\it There exists a positive constant $M_{1}$, independent of 
$\varepsilon \in (0, 1]$, such that }
\begin{equation}
\frac{1}{2}\int_{0}^{T}\bigl|\boldsymbol{v}_{\varepsilon }(s)\bigr|_{\boldsymbol{V}_{0}}^{2}ds
+\int_{0}^{T}\int_{\Omega }\widehat{\beta }_{\varepsilon }\bigl(v_{\varepsilon }(s)+m_{0}\bigr)dxds
+\int_{0}^{T}\int_{\Gamma }\widehat{\beta }_{\Gamma , \varepsilon }\bigl(v_{\Gamma , \varepsilon }(s)+m_{0}\bigr)d\Gamma ds
\leq M_{1}. \label{es1}
\end{equation}
\paragraph{Proof}
From (\ref{pi}), (\ref{pig}) and the assumption (A3), 
note that $\widetilde{\pi }, \widetilde{\pi }_{\Gamma }$ is 
grobally {L}ipschitz continuous on $\mathbb{R}$. 
We denote the {L}ipschitz constant of $\widetilde{\pi }, \widetilde{\pi }_{\Gamma }$ by $\widetilde{L}, \widetilde{L}_{\Gamma }$, respectively. 
Moreover, we can take primitive functions $\widehat{\widetilde{\pi }}$ 
and $\widehat{\widetilde{\pi }}_{\Gamma }$ of $\widetilde{\pi }$ and $\widetilde{\pi }_{\Gamma }$ 
satisfying 
$$\int_{\Omega }\widehat{\widetilde{\pi }}\bigl(v_{\varepsilon }(s)\bigr)dx \geq 0, 
\quad \int_{\Gamma }\widehat{\widetilde{\pi }}_{\Gamma }(v_{\Gamma , \varepsilon }(s))d\Gamma \geq 0$$
for a.a.\ $s\in (0, T)$, 
respectively. 
Now, we test (\ref{apo1}) at time $s\in (0, T)$ by $\boldsymbol{v}_{\varepsilon }(s)$ 
and use the {Y}oung inequality. 
Then, we deduce that 
\begin{eqnarray*}
&&\frac{1}{2}\frac{d}{ds}\bigl|\boldsymbol{v}_{\varepsilon }(s)\bigr|_{\boldsymbol{V}_{0}^{*}}^{2}
+\frac{\varepsilon }{2}\frac{d}{ds}\bigl|\boldsymbol{v}_{\varepsilon }(s)\bigr|_{\boldsymbol{H}_{0}}^{2}
+\varphi _{\varepsilon }\bigl(\boldsymbol{v}_{\varepsilon }(s)\bigr) \nonumber \\
&&\quad \leq \bigl(\boldsymbol{P}\boldsymbol{f}(s)-\boldsymbol{P}\bigl(\widetilde{\boldsymbol{\pi }}
\bigl(\boldsymbol{v}_{\varepsilon }(s)+m_{0}\boldsymbol{1}\bigr)\bigr), \boldsymbol{v}_{\varepsilon }(s)\bigr)_{\boldsymbol{H}_{0}}
+\varphi _{\varepsilon }(0) \nonumber \\
&&\quad \leq c_{{\rm P}}\bigl|\boldsymbol{f}(s)\bigr|_{\boldsymbol{H}}^{2}
+\frac{1}{4c_{{\rm P}}}\bigl|\boldsymbol{v}_{\varepsilon }(s)\bigr|_{\boldsymbol{H}_{0}}^{2}
+c_{{\rm P}}M+\varphi (0) \nonumber \\
&&\quad \leq \frac{1}{4}\bigl|\boldsymbol{v}_{\varepsilon }(s)\bigr|_{\boldsymbol{V}_{0}}^{2}
+c_{{\rm P}}\bigl|\boldsymbol{f}(s)\bigr|_{\boldsymbol{H}}^{2}
+c_{{\rm P}}M+\varphi (0) \nonumber \\
&&\quad \leq \frac{1}{2}\varphi _{\varepsilon }\bigl(\boldsymbol{v}_{\varepsilon }(s)\bigr)
+c_{{\rm P}}\bigl|\boldsymbol{f}(s)\bigr|_{\boldsymbol{H}}^{2}
+c_{{\rm P}}M+\varphi (0) \label{ine3}
\end{eqnarray*}
for a.a.\ $s\in (0, T)$. 
Namely, we have 
\begin{eqnarray*}
&&\frac{1}{2}\frac{d}{ds}\bigl|\boldsymbol{v}_{\varepsilon }(s)\bigr|_{\boldsymbol{V}_{0}^{*}}^{2}
+\frac{\varepsilon }{2}\frac{d}{ds}\bigl|\boldsymbol{v}_{\varepsilon }(s)\bigr|_{\boldsymbol{H}_{0}}^{2}
+\frac{1}{4}\bigl|\boldsymbol{v}_{\varepsilon }(s)\bigr|_{\boldsymbol{V}_{0}}^{2} \nonumber \\
&&\quad \quad \quad +\frac{1}{2}\int_{\Omega }\widehat{\beta }_{\varepsilon }\bigl(v_{\varepsilon }(s)+m_{0}\bigr)dx
+\frac{1}{2}\int_{\Gamma }\widehat{\beta }_{\Gamma , \varepsilon }\bigl(v_{\Gamma , \varepsilon }(s)+m_{0}\bigr)d\Gamma \nonumber \\
&&\quad \leq c_{{\rm P}}\bigl|\boldsymbol{f}(s)\bigr|_{\boldsymbol{H}}^{2}
+c_{{\rm P}}M+\varphi (0) \label{mulv}
\end{eqnarray*}
for a.a.\ $s\in (0, T)$. 
Integrating it over $(0, T)$ and using the periodic property, we see that 
\begin{eqnarray*}
&&\frac{1}{2}\int_{0}^{T}\bigl|\boldsymbol{v}_{\varepsilon }(s)\bigr|_{\boldsymbol{V}_{0}}^{2}
+\int_{0}^{T}\int_{\Omega }\widehat{\beta }_{\varepsilon }\bigl(v_{\varepsilon }(s)+m_{0}\bigr)dxds
+\int_{0}^{T}\int_{\Gamma }\widehat{\beta }_{\Gamma , \varepsilon }\bigl(v_{\Gamma , \varepsilon }(s)+m_{0}\bigr)d\Gamma ds \nonumber \\
&&\quad \leq 2c_{{\rm P}}\bigl|\boldsymbol{f}\bigr|_{L^{2}(0, T; \boldsymbol{H})}^{2}
+2c_{{\rm P}}TM+2T\varphi (0). \label{m2}
\end{eqnarray*}
This yields that the estimete (\ref{es1}) holds. 
\hfill $\Box $

%%%second estimate%%%%%%%%%%%%%%%%%%%%%%%%%%%%%%%%%%%%%%%%%%
\paragraph{Lemma 3.3.}
{\it There exists a positive constant $M_{2}$, independent of 
$\varepsilon \in (0, 1]$, such that }
\begin{equation*}
\varepsilon \int_{0}^{T}\bigl|\boldsymbol{v}_{\varepsilon }'(s)\bigr|_{\boldsymbol{H}_{0}}^{2}ds
+\frac{1}{2}\int_{0}^{T}\bigl|\boldsymbol{v}_{\varepsilon }'(s)\bigr|_{\boldsymbol{V}_{0}^{*}}^{2}ds
\leq M_{2}. \label{es2}
\end{equation*}
\paragraph{Proof}
We test (\ref{apo1}) at time $s\in (0, T)$ by $\boldsymbol{v}_{\varepsilon }'(s)$. 
Then, by using the {Y}oung inequality, 
we see that 
\begin{eqnarray*}
&&\varepsilon \bigl|\boldsymbol{v}_{\varepsilon }'(s)\bigr|_{\boldsymbol{H}_{0}}^{2}+
\bigl|\boldsymbol{v}_{\varepsilon }'(s)\bigr|_{\boldsymbol{V}_{0}^{*}}^{2}
+\frac{d}{ds}\varphi _{\varepsilon }\bigl(\boldsymbol{v}_{\varepsilon }(s)\bigr) \\
&&\quad \quad \quad +\frac{d}{ds}\int_{\Omega }\widehat{\widetilde{\pi }}\bigl(v_{\varepsilon }(s)+m_{0}\bigr)dx
+\frac{d}{ds}\int_{\Gamma }\widehat{\widetilde{\pi }}_{\Gamma }\bigl(v_{\Gamma , \varepsilon }(s)+m_{0}\bigr)d\Gamma  \\
&&\quad =\bigl(\boldsymbol{Pf}(s), \boldsymbol{v}_{\varepsilon }'(s)\bigr)_{\boldsymbol{H}_{0}} \\
&&\quad \leq \frac{1}{2}\bigl|\boldsymbol{f}(s)\bigr|_{\boldsymbol{V}}^{2}
+\frac{1}{2}\bigl|\boldsymbol{v}_{\varepsilon }'(s)\bigr|_{\boldsymbol{V}_{0}^{*}}^{2}
\end{eqnarray*}
for a.a.\ $s\in (0, T)$. 
This implies that 
\begin{eqnarray}
&&\varepsilon \bigl|\boldsymbol{v}_{\varepsilon }'(s)\bigr|_{\boldsymbol{H}_{0}}^{2}
+\frac{1}{2}\bigl|\boldsymbol{v}_{\varepsilon }'(s)\bigr|_{\boldsymbol{V}_{0}^{*}}^{2}
+\frac{d}{ds}\varphi _{\varepsilon }\bigl(\boldsymbol{v}_{\varepsilon }(s)\bigr) \nonumber \\
&&\quad \quad \quad +\frac{d}{ds}\int_{\Omega }\widehat{\widetilde{\pi }}\bigl(v_{\varepsilon }(s)+m_{0}\bigr)dx
+\frac{d}{ds}\int_{\Gamma }\widehat{\widetilde{\pi }}_{\Gamma }\bigl(v_{\Gamma , \varepsilon }(s)+m_{0}\bigr)d\Gamma \nonumber \\
&&\quad \leq \frac{1}{2}\bigl|\boldsymbol{f}(s)\bigr|_{\boldsymbol{V}}^{2} \label{mulvd}
\end{eqnarray}
for a.a.\ $s\in (0, T)$. 
Therefore, by integrating it over $(0, T)$ with respect to $s$ and using the periodic property, 
we can conclude. 
\hfill $\Box $

%%%third estimate%%%%%%%%%%%%%%%%%%%%%%%%%%%%%%%%%%%%%%%%%%
\paragraph{Lemma 3.4.}
{\it There exists a positive constant $M_{3}$, independent of 
$\varepsilon \in (0, 1]$, such that }
\begin{equation}\label{es3}
\frac{1}{2}\bigl|\boldsymbol{v}_{\varepsilon }(t)\bigr|_{\boldsymbol{V}_{0}}^{2}
+\int_{\Omega }\widehat{\widetilde{\pi }}\bigl(v_{\varepsilon }(t)+m_{0}\bigr)dx
+\int_{\Gamma }\widehat{\widetilde{\pi }}_{\Gamma }\bigl(v_{\Gamma , \varepsilon }(t)+m_{0}\bigr)d\Gamma 
\leq M_{3} 
\end{equation}
{\it for a.a.\ $t\in [0, T]$. }
\paragraph{Proof}
For each $s, t\in [0, T]$ such that $s\leq t$, 
we integrate (\ref{mulvd}) over $[s, t]$. 
Then, by neglecting the first two positive terms, we see that 
\begin{eqnarray*}
&&\varphi _{\varepsilon }\bigl(\boldsymbol{v}_{\varepsilon }(t)\bigr)
+\int_{\Omega }\widehat{\widetilde{\pi }}\bigl(v_{\varepsilon }(t)+m_{0}\bigr)dx
+\int_{\Gamma }\widehat{\widetilde{\pi }}_{\Gamma }\bigl(v_{\Gamma , \varepsilon }(t)+m_{0}\bigr)d\Gamma \\
&&\quad \leq \varphi _{\varepsilon }\bigl(\boldsymbol{v}_{\varepsilon }(s)\bigr)
+\int_{\Omega }\widehat{\widetilde{\pi }}\bigl(v_{\varepsilon }(s)+m_{0}\bigr)dx
+\int_{\Gamma }\widehat{\widetilde{\pi }}_{\Gamma }\bigl(v_{\Gamma , \varepsilon }(s)+m_{0}\bigr)d\Gamma 
+\frac{1}{2}\int_{0}^{T}\bigl|\boldsymbol{f}(s)\bigr|_{\boldsymbol{V}_{0}}^{2}ds
\end{eqnarray*}
for a.a.\ $s, t\in [0, T]$. 
Now, Integrating it over $(0, t)$ with respect to $s$, it follows that 
\begin{eqnarray}
&&\frac{t}{2}\bigl|\boldsymbol{v}_{\varepsilon }(t)\bigr|_{\boldsymbol{V}_{0}}^{2}
+t\int_{\Omega }\widehat{\beta }_{\varepsilon }\bigl(v_{\varepsilon }(t)+m_{0}\bigr)dx
+t\int_{\Gamma }\widehat{\beta }_{\Gamma , \varepsilon }\bigl(v_{\Gamma , \varepsilon }(t)+m_{0}\bigr)d\Gamma \nonumber \\
&&\quad \leq \frac{1}{2}\int_{0}^{T}\bigl|\boldsymbol{v}_{\varepsilon }(s)\bigr|_{\boldsymbol{V}_{0}}^{2}ds
+\int_{0}^{T}\int_{\Omega }\widehat{\beta }_{\varepsilon }\bigl(v_{\varepsilon }(s)+m_{0}\bigr)dxds
+\int_{0}^{T}\int_{\Gamma }\widehat{\beta }_{\Gamma , \varepsilon }\bigl(v_{\Gamma , \varepsilon }(s)+m_{0}\bigr)d\Gamma ds \nonumber \\
&&\quad \quad \quad +\int_{0}^{T}\int_{\Omega }\widehat{\widetilde{\pi }}\bigl(v_{\varepsilon }(s)+m_{0}\bigr)dxds
+\int_{0}^{T}\int_{\Gamma }\widehat{\widetilde{\pi }}_{\Gamma }\bigl(v_{\Gamma , \varepsilon }(s)+m_{0}\bigr)d\Gamma ds \nonumber \\
&&\quad \quad \quad \quad \quad +\frac{T}{2}\int_{0}^{T}\bigl|\boldsymbol{f}(s)\bigr|_{\boldsymbol{V}_{0}}^{2}ds \label{ints}
\end{eqnarray}
for a.a.\ $t\in [0, T]$. 
Here, Note that we have 
\begin{eqnarray*}
\bigl|\widehat{\widetilde{\pi }}(r)\bigr|&\leq &\int_{0}^{r}\bigl|\widetilde{\pi }(\tau )\bigr|d\tau \\
                                         &\leq &\widetilde{L}\int_{0}^{r}|\tau |d\tau +\int_{0}^{r}\bigl|\widetilde{\pi }(0)\bigr|d\tau \\
                                         &=&\frac{\widetilde{L}}{2}r^{2}+\bigl|\widetilde{\pi }(0)\bigr|r %\quad {\rm for~all~}r\in \mathbb{R}, 
\end{eqnarray*}
for all $r\in \mathbb{R}$. 
and similarly, 
\begin{eqnarray*}
\bigl|\widehat{\widetilde{\pi }}_{\Gamma }(r)\bigr|\leq \frac{\widetilde{L}_{\Gamma }}{2}r^{2}+\bigl|\widetilde{\pi }_{\Gamma }(0)\bigr|r \quad {\rm for~all~}r\in \mathbb{R}. 
\end{eqnarray*}
Then, by using the {Y}oung inequality, we infer that 
\begin{eqnarray}
\int_{\Omega }\widehat{\widetilde{\pi }}\bigl(v_{\varepsilon }(s)+m_{0}\bigr)dx
&\leq &\int_{\Omega }\left(\frac{\widetilde{L}}{2}\bigl|v_{\varepsilon }(s)+m_{0}\bigr|^{2}+
\bigl|\widetilde{\pi }(0)\bigr|\bigl|v_{\varepsilon }(s)+m_{0}\bigr|\right)dx \nonumber \\
&\leq &\widetilde{L}\int_{\Omega }\bigl|v_{\varepsilon }(s)+m_{0}\bigr|^{2}dx
+\frac{1}{2\widetilde{L}}\bigl|\widetilde{\pi }(0)\bigr|^{2}|\Omega | \nonumber \\
&\leq &2\widetilde{L}\int_{\Omega }\bigl|v_{\varepsilon }(s)\bigr|^{2}dx
+2m_{0}^{2}|\Omega |+\frac{1}{2\widetilde{L}}\bigl|\widetilde{\pi }(0)\bigr|^{2}|\Omega | \label{pihat}
\end{eqnarray}
for a.a.\ $s\in [0, T]$. 
Similarly, we have 
\begin{eqnarray}
\int_{\Gamma }\widehat{\widetilde{\pi }}_{\Gamma }\bigl(v_{\Gamma , \varepsilon }(s)+m_{0}\bigr)d\Gamma 
\leq 2\widetilde{L}_{\Gamma }\int_{\Gamma }\bigl|v_{\Gamma , \varepsilon }(s)\bigr|^{2}d\Gamma 
+2m_{0}^{2}|\Gamma |+\frac{1}{2\widetilde{L}_{\Gamma }}\bigl|\widetilde{\pi }_{\Gamma }(0)\bigr|^{2}|\Gamma | \label{pighat}
\end{eqnarray}
for a.a.\ $s\in [0, T]$. 
Thus, on account of (\ref{ints})--(\ref{pighat}), we deduce that 
\begin{eqnarray*}
&&\frac{t}{2}\bigl|\boldsymbol{v}_{\varepsilon }(t)\bigr|_{\boldsymbol{V}_{0}}^{2}
+t\int_{\Omega }\widehat{\beta }_{\varepsilon }\bigl(v_{\varepsilon }(t)+m_{0}\bigr)dx
+t\int_{\Gamma }\widehat{\beta }_{\Gamma , \varepsilon }\bigl(v_{\Gamma , \varepsilon }(t)+m_{0}\bigr)d\Gamma \\
&&\quad \leq \frac{1}{2}\int_{0}^{T}\bigl|\boldsymbol{v}_{\varepsilon }(s)\bigr|_{\boldsymbol{V}_{0}}^{2}ds
+\int_{0}^{T}\int_{\Omega }\widehat{\beta }_{\varepsilon }\bigl(v_{\varepsilon }(s)+m_{0}\bigr)dxds
+\int_{0}^{T}\int_{\Gamma }\widehat{\beta }_{\Gamma , \varepsilon }\bigl(v_{\Gamma , \varepsilon }(s)+m_{0}\bigr)d\Gamma ds \\
&&\quad \quad \quad +2\widetilde{L}\int_{\Omega }\bigl|v_{\varepsilon }(s)\bigr|^{2}dx
+2\widetilde{L}_{\Gamma }\int_{\Gamma }\bigl|v_{\Gamma , \varepsilon }(s)\bigr|^{2}d\Gamma 
+\frac{T}{2}\int_{0}^{T}\bigl|\boldsymbol{f}(s)\bigr|_{\boldsymbol{V}_{0}}^{2}ds
+\tilde{M}_{4} \\
&&\quad \leq \left(\frac{1}{2}+\widehat{L}c_{{\rm P}}\right)\int_{0}^{T}\bigl|\boldsymbol{v}_{\varepsilon }(s)\bigr|_{\boldsymbol{V}_{0}}^{2}ds 
+\int_{0}^{T}\int_{\Omega }\widehat{\beta }_{\varepsilon }\bigl(v_{\varepsilon }(s)+m_{0}\bigr)dxds \\
&&\quad \quad \quad +\int_{0}^{T}\int_{\Gamma }\widehat{\beta }_{\Gamma , \varepsilon }\bigl(v_{\Gamma , \varepsilon }(s)+m_{0}\bigr)d\Gamma ds 
+\frac{T}{2}\int_{0}^{T}\bigl|\boldsymbol{f}(s)\bigr|_{\boldsymbol{V}_{0}}^{2}ds
+\tilde{M}_{4}
\end{eqnarray*}
for a.a.\ $t\in [0, T]$, where $\widehat{L}:=\max \{2\widetilde{L}, 2\widetilde{L}_{\Gamma }\}$ and 
\begin{equation*}
\tilde{M}_{4}:=2m_{0}^{2}|\Omega |+\frac{1}{2\widetilde{L}}\bigl|\widetilde{\pi }(0)\bigr|^{2}|\Omega |
+2m_{0}^{2}|\Gamma |+\frac{1}{2\widetilde{L}_{\Gamma }}\bigl|\widetilde{\pi }_{\Gamma }(0)\bigr|^{2}|\Gamma |. 
\end{equation*}
In particular, putting $t:=T$ and dividing it by $T$, it follows that 
\begin{eqnarray}
&&\frac{1}{2}\bigl|\boldsymbol{v}_{\varepsilon }(T)\bigr|_{\boldsymbol{V}_{0}}^{2}
+\int_{\Omega }\widehat{\beta }_{\varepsilon }\bigl(v_{\varepsilon }(T)+m_{0}\bigr)dx
+\int_{\Gamma }\widehat{\beta }_{\Gamma , \varepsilon }\bigl(v_{\Gamma , \varepsilon }(T)+m_{0}\bigr)d\Gamma \nonumber \\
&&\quad \leq \frac{1}{T}\left(\frac{1}{2}+\widehat{L}c_{{\rm P}}\right)\int_{0}^{T}\bigl|\boldsymbol{v}_{\varepsilon }(s)\bigr|_{\boldsymbol{V}_{0}}^{2}ds 
+\frac{1}{T}\int_{0}^{T}\int_{\Omega }\widehat{\beta }_{\varepsilon }\bigl(v_{\varepsilon }(s)+m_{0}\bigr)dxds \nonumber \\
&&\quad \quad \quad +\frac{1}{T}\int_{0}^{T}\int_{\Gamma }\widehat{\beta }_{\Gamma , \varepsilon }\bigl(v_{\Gamma , \varepsilon }(s)+m_{0}\bigr)d\Gamma ds 
+\frac{1}{2}\int_{0}^{T}\bigl|\boldsymbol{f}(s)\bigr|_{\boldsymbol{V}_{0}}^{2}ds
+\frac{\tilde{M}_{4}}{T}.  \label{estiT}
\end{eqnarray}
Combining (\ref{es1}) and (\ref{estiT}), there exists a positive constant $\tilde{M}_{3}$ such that 
\begin{equation*}
\frac{1}{2}\bigl|\boldsymbol{v}_{\varepsilon }(T)\bigr|_{\boldsymbol{V}_{0}}^{2}
+\int_{\Omega }\widehat{\beta }_{\varepsilon }\bigl(v_{\varepsilon }(T)+m_{0}\bigr)dx
+\int_{\Gamma }\widehat{\beta }_{\Gamma , \varepsilon }\bigl(v_{\Gamma , \varepsilon }(T)+m_{0}\bigr)d\Gamma 
\leq \tilde{M}_{3}. 
\end{equation*}
From the periodic property, we have 
\begin{equation}
\varphi _{\varepsilon }\bigl(\boldsymbol{v}_{\varepsilon }(0)\bigr)=
\frac{1}{2}\bigl|\boldsymbol{v}_{\varepsilon }(0)\bigr|_{\boldsymbol{V}_{0}}^{2}
+\int_{\Omega }\widehat{\beta }_{\varepsilon }\bigl(v_{\varepsilon }(0)+m_{0}\bigr)dx
+\int_{\Gamma }\widehat{\beta }_{\Gamma , \varepsilon }\bigl(v_{\Gamma , \varepsilon }(0)+m_{0}\bigr)d\Gamma 
\leq \tilde{M}_{3}. \label{tilm3}
\end{equation}
Now, integrating (\ref{mulvd}) by $(0, t)$ with respect to $s$, 
it follows from (\ref{pihat})--(\ref{pighat}) that 
\begin{eqnarray}
&&\varphi _{\varepsilon }\bigl(\boldsymbol{v}_{\varepsilon }(t)\bigr)
+\int_{\Omega }\widehat{\widetilde{\pi }}\bigl(v_{\varepsilon }(t)+m_{0}\bigr)dx
+\int_{\Gamma }\widehat{\widetilde{\pi }}_{\Gamma }\bigl(v_{\Gamma , \varepsilon }(t)+m_{0}\bigr)d\Gamma \nonumber \\
&&\quad \leq \varphi _{\varepsilon }\bigl(\boldsymbol{v}_{\varepsilon }(0)\bigr)
+\int_{\Omega }\widehat{\widetilde{\pi }}\bigl(v_{\varepsilon }(0)+m_{0}\bigr)dx
+\int_{\Gamma }\widehat{\widetilde{\pi }}_{\Gamma }\bigl(v_{\Gamma , \varepsilon }(0)+m_{0}\bigr)d\Gamma 
+\frac{1}{2}\int_{0}^{T}\bigl|\boldsymbol{f}(s)\bigr|_{\boldsymbol{V}_{0}}^{2}ds \nonumber \\
&&\quad \leq \bigl(1+2\widehat{L}c_{{\rm P}}\bigr)\varphi _{\varepsilon }\bigl(\boldsymbol{v}_{\varepsilon }(0)\bigr)
+\frac{1}{2}\int_{0}^{T}\bigl|\boldsymbol{f}(s)\bigr|_{\boldsymbol{V}_{0}}^{2}ds
+\tilde{M}_{4} \label{esv0}
\end{eqnarray}
for a.a.\ $t\in [0, T]$. 
Therefore, by virtue of (\ref{tilm3})--(\ref{esv0}), there exists a positive constnat $M_{3}$ such that 
the estimate (\ref{es3}) holds. 
\hfill $\Box $

%%%forth estimate%%%%%%%%%%%%%%%%%%%%%%%%%%%%%%%%%%%%%%%%%%
\paragraph{Lemma 3.5.}
{\it There exists a positive constant $M_{4}$, independent of 
$\varepsilon \in (0, 1]$, such that }
\begin{equation}\label{es4}
\delta _{0}\int_{0}^{T}\bigl|\beta _{\varepsilon }\bigl(v_{\varepsilon }(s)+m_{0}\bigr)\bigr|_{L^{1}(\Omega )}^{2}ds
+\delta _{0}\int_{0}^{T}\bigl|\beta _{\Gamma , \varepsilon }\bigl(v_{\Gamma , \varepsilon }(s)+m_{0}\bigr)\bigr|_{L^{1}(\Gamma )}^{2}ds
\leq M_{4} 
\end{equation}
for some positive constants $\delta _{0}$. 

\paragraph{Proof}
To show this lemma, 
we can employ the method of \cite[Lemma 4.1, 4.3]{CF15}, 
because of imposing same assumptions as \cite{CF15} for $\beta , \beta _{\Gamma }$ and being $m_{0}\in {\rm int}D(\beta _{\Gamma })$. 
Therefore, we can also exploit the following inequalities stated in \cite[Sect.\ 5]{GMS09}: 
for each $\varepsilon \in (0, 1]$, 
there exist two positive constants $\delta _{0}$ and $c_{1}$ such that 
$$\beta _{\varepsilon }(r)(r-m_{0})\geq \delta _{0}\bigl|\beta _{\varepsilon }(r)\bigr|-c_{1}, 
\quad \beta _{\Gamma , \varepsilon }(r)(r-m_{0})\geq \delta _{0}\bigl|\beta _{\Gamma , \varepsilon }(r)\bigr|-c_{1}$$
for all $r\in \mathbb{R}$. 
Hence, it follows that 
\begin{equation}\label{gms}
\bigl(\boldsymbol{\beta }_{\varepsilon }\bigl(\boldsymbol{u}_{\varepsilon }(s)\bigr), \boldsymbol{v}_{\varepsilon }(s)\bigr)_{\boldsymbol{H}}
\geq \delta _{0}\int_{\Omega }\bigl|\beta _{\varepsilon }\bigl(u_{\varepsilon }(s)\bigr)\bigr|dx-c_{1}|\Omega |
+\delta _{0}\int_{\Gamma }\bigl|\beta _{\Gamma , \varepsilon }\bigl(u_{\Gamma , \varepsilon }(s)\bigr)\bigr|d\Gamma -c_{1}|\Gamma |
\end{equation}
for a.a.\ $s\in (0, T)$. 
On the other hand, 
we test (\ref{apo1}) at time $s\in (0, T)$ by $\boldsymbol{v}_{\varepsilon }(s)$. 
Then, from (\ref{subvarp}), we see that 
\begin{eqnarray}
&&\bigl(\varepsilon \boldsymbol{v}_{\varepsilon }'(s), \boldsymbol{v}_{\varepsilon }(s)\bigr)_{\boldsymbol{H}_{0}}
+\bigl(\boldsymbol{v}_{\varepsilon }'(s), \boldsymbol{v}_{\varepsilon }(s)\bigr)_{\boldsymbol{V}_{0}^{*}}
+\bigl(\boldsymbol{A}\boldsymbol{v}_{\varepsilon }(s), \boldsymbol{v}_{\varepsilon }(s)\bigr)_{\boldsymbol{H}_{0}}
+\bigl(\boldsymbol{P}\boldsymbol{\beta }_{\varepsilon }\bigl(\boldsymbol{u}_{\varepsilon }(s)\bigr), \boldsymbol{v}_{\varepsilon }(s)\bigr)_{\boldsymbol{H}_{0}} \nonumber \\
&&\quad \leq \bigl(\boldsymbol{f}(s)-\widetilde{\boldsymbol{\pi }}\bigl(\boldsymbol{u}_{\varepsilon }(s)\bigr), \boldsymbol{v}_{\varepsilon }(s)\bigr)_{\boldsymbol{H}}. \label{inev}
\end{eqnarray}
Hence, from (\ref{gms})--(\ref{inev}) and the maximal monotonicity of $\boldsymbol{A}$, by squaring we have 
\begin{eqnarray*}
&&\left(\delta _{0}\int_{\Omega }\bigl|\beta _{\varepsilon }\bigl(u_{\varepsilon }(s)\bigr)\bigr|dx+
\delta _{0}\int_{\Gamma }\bigl|\beta _{\Gamma , \varepsilon }\bigl(u_{\Gamma , \varepsilon }(s)\bigr)\bigr|d\Gamma \right)^{2}
\leq 3c_{1}^{2}(|\Omega |+|\Gamma |)^{2} \nonumber \\
&&\quad +9\bigl(\bigl|\boldsymbol{f}(s)\bigr|_{\boldsymbol{H}}^{2}+\bigl|\boldsymbol{\pi }\bigl(\boldsymbol{u}_{\varepsilon }(s)\bigr)\bigr|_{\boldsymbol{H}}^{2}
+\varepsilon ^{2}\bigl|\boldsymbol{v}_{\varepsilon }'(s)\bigr|_{\boldsymbol{H}_{0}}^{2}\bigr)\bigl|\boldsymbol{v}_{\varepsilon }(s)\bigr|_{\boldsymbol{H}_{0}}^{2}
+3\bigl|\boldsymbol{v}_{\varepsilon }'(s)\bigr|_{\boldsymbol{V}_{0}^{*}}^{2}\bigl|\boldsymbol{v}_{\varepsilon }(s)\bigr|_{\boldsymbol{V}_{0}^{*}}^{2}
\end{eqnarray*}
for a.a.\ $s\in (0, T)$. 
Therefore, from the {L}ipschitz continuity of $\widetilde{\pi }, \widetilde{\pi }_{\Gamma }$ and Lemma 3.4, 
by integrating it over $(0, T)$ with respect to $s$, 
there exists a positive constant $M_{4}$ such that 
the estimate (\ref{es4}) holds. 
\hfill $\Box $

%Following lemmas are proved exactly the same as in \cite[Lemmas 4.3, 4.4 and 4.5]{CF15}
%because the form $\boldsymbol{\mu }_{\varepsilon }$ and $\omega _{\varepsilon }$ are the same as \cite{CF15}
%and we already obtained estimates respect to $\boldsymbol{v}_{\varepsilon }$ and $\boldsymbol{\beta }_{\varepsilon }$ enough to prove. 
%The key point to prove lemmas is that we can adapt the theory of the elliptic regularity (see, e.g.\, \cite[Theorem 3.2, p.\ 1.79]{BG87})
%becuse we have the term $\Delta _{\Gamma }v _{\varepsilon }$ in the approximate sequence $\mu _{\Gamma , \varepsilon }$. 
%This is the big advantage of dynamic boundary condition (\ref{eq5}). 
%So, we omit the proof the following lemmas. 

%%%fifth estimate%%%%%%%%%%%%%%%%%%%%%%%%%%%%%%%%%%%%%%%%%%
\paragraph{Lemma 3.6.}
{\it There exists a positive constants $M_{5}$, independent of 
$\varepsilon \in (0, 1]$, such that }
\begin{gather}\label{es5}
\int_{0}^{T}\bigl|\boldsymbol{\mu }_{\varepsilon }(s)\bigr|_{\boldsymbol{V}}^{2}ds\leq M_{5}. 
%\int_{0}^{T}\bigl|\boldsymbol{\beta }_{\varepsilon }\bigl(u_{\varepsilon }(s)\bigr)\bigr|_{H}^{2}ds
%+\int_{0}^{T}\bigl|\boldsymbol{\beta }_{\Gamma , \varepsilon }\bigl(u_{\Gamma , \varepsilon }(s)\bigr)\bigr|_{H_{\Gamma }}^{2}ds\leq M_{6}, \\
%\int_{0}^{T}\bigl|\boldsymbol{v}_{\varepsilon }(s)\bigr|_{\boldsymbol{W}}^{2}ds\leq M_{7}
\end{gather}
%{\it for a.a.\ $s\in (0, T)$. }
\paragraph{Proof}
Firstly, %we deduce the uniform estimate of the volume of $\widetilde{\boldsymbol{\pi }}(\boldsymbol{v}_{\varepsilon }(s)+m_{0}\boldsymbol{1})$. 
by using the {L}ipschitz continuity of $\widetilde{\pi }, \widetilde{\pi }_{\Gamma }$ and the H\"older inequality, 
it follows from (\ref{PW0}) and Lemma 3.4 that 
there exists a positive constnat $M_{5}^{*}$ such that 
\begin{eqnarray}\label{volpi}
&&\bigl|m\bigl(\widetilde{\boldsymbol{\pi }}\bigl(\boldsymbol{v}_{\varepsilon }(s)+m_{0}\boldsymbol{1}\bigr)\bigr)\bigr| \nonumber \\
&&\quad \leq \frac{1}{|\Omega |+|\Gamma |}\left\{\int_{\Omega }\bigl|\widetilde{\pi }\bigl(v_{\varepsilon }(s)+m_{0}\bigr)\bigr|dx
+\int_{\Gamma }\bigl|\widehat{\pi }_{\Gamma }\bigl(v_{\Gamma , \varepsilon }(s)+m_{0}\bigr)\bigr|d\Gamma \right\} \nonumber \\
&&\quad \leq \frac{1}{|\Omega |+|\Gamma |}\left\{\widetilde{L}|\Omega |^{\frac{1}{2}}\bigl|v_{\varepsilon }(s)\bigr|_{H}^{2}
+\widetilde{L}|\Omega ||m_{0}|+|\Omega |\bigl|\widetilde{\pi }(0)\bigr| \right. \nonumber \\
&&\left. \quad \quad \quad +\widetilde{L}_{\Gamma }|\Gamma |^{\frac{1}{2}}\bigl|v_{\Gamma , \varepsilon }(s)\bigr|_{H_{\Gamma }}^{2}
+\widetilde{L}_{\Gamma }|\Gamma ||m_{0}|+|\Gamma |\bigl|\widetilde{\pi }_{\Gamma }(0)\bigr|\right\} \nonumber \\
&&\quad \leq \frac{1}{|\Omega |+|\Gamma |}M_{5}^{*}\left\{\bigl|\boldsymbol{v}_{\varepsilon }(s)\bigr|_{\boldsymbol{V}_{0}}^{2}+1\right\} \nonumber \\
&&\quad \leq \frac{1}{|\Omega |+|\Gamma |}M_{5}^{*}(M_{3}+1)=:\tilde{M}_{5} 
\end{eqnarray}
for a.a.\ $s\in (0, T)$. 
Therefore, owing to (\ref{volpi}) we deduce that 
\begin{eqnarray*}
\bigl|m\bigl(\boldsymbol{\mu }_{\varepsilon }(s)\bigr)\bigr|^{2}
&=&\bigl|m\bigl(\widetilde{\boldsymbol{\pi }}\bigl(\boldsymbol{v}_{\varepsilon }(s)+m_{0}\boldsymbol{1}\bigr)-\boldsymbol{f}(s)\bigr)\bigr|^{2} \\
&\leq &2\tilde{M}_{5}^{2}+\frac{4}{(|\Omega |+|\Gamma |)^{2}}\bigl(\bigl|f(s)\bigr|_{L^{1}(\Omega )}+\bigl|f_{\Gamma }(s)\bigr|_{L^{1}(\Gamma )}\bigr)=:\hat{M}_{5}
\end{eqnarray*}
for a.a.\ $s\in (0, T)$. 
Next, from (\ref{PW0}), (\ref{vmu}) and the fact $\boldsymbol{P}\boldsymbol{\mu }_{\varepsilon }(s)=\boldsymbol{\mu }_{\varepsilon }(s)-m(\boldsymbol{\mu }_{\varepsilon }(s))\boldsymbol{1}$ 
for a.a.\ $s\in (0, T)$, we deduce that 
\begin{eqnarray*}
\int_{0}^{T}\bigl|\boldsymbol{\mu }_{\varepsilon }(s)\bigr|_{\boldsymbol{V}}^{2}ds
&\leq &2\int_{0}^{T}\bigl|\boldsymbol{P}\boldsymbol{\mu }_{\varepsilon }(s)\bigr|_{\boldsymbol{V}}^{2}ds
+2\int_{0}^{T}\bigl|m\bigl(\boldsymbol{\mu }_{\varepsilon }(s)\bigr)\boldsymbol{1}\bigr|_{\boldsymbol{V}}^{2}ds \\
&\leq &2c_{{\rm P}}\int_{0}^{T}\bigl|\boldsymbol{P}\boldsymbol{\mu }_{\varepsilon }(s)\bigr|_{\boldsymbol{V}_{0}}^{2}ds
+2(|\Omega |+|\Gamma |)\int_{0}^{T}\bigl|m\bigl(\boldsymbol{\mu }_{\varepsilon }(s)\bigr)\bigr|^{2}ds \\
&\leq &2c_{{\rm P}}\int_{0}^{T}\bigl|\boldsymbol{v}_{\varepsilon }'(s)\bigr|_{\boldsymbol{V}_{0}^{*}}^{2}ds
+2T(|\Omega |+|\Gamma |)\hat{M}_{5}^{2}. \\
\end{eqnarray*}
Thus, from Lemma 3.3, there exist a positive constant $M_{5}$ such that 
the estimate (\ref{es5}) holds. 
\hfill $\Box $

%%%sixth estimate%%%%%%%%%%%%%%%%%%%%%%%%%%%%%%%%%%%%%%%%%%
\paragraph{Lemma 3.7.}
{\it There exists a positive constant $M_{6}$, independent of 
$\varepsilon \in (0, 1]$, such that }
\begin{gather}
%\int_{0}^{T}\bigl|\boldsymbol{\mu }_{\varepsilon }(s)\bigr|_{\boldsymbol{V}}^{2}ds\leq M_{5}
\frac{1}{2}\int_{0}^{T}\bigl|\beta _{\varepsilon }\bigl(v_{\varepsilon }(s)+m_{0}\bigr)\bigr|_{H}^{2}ds
+\frac{1}{4\rho }\int_{0}^{T}\bigl|\beta _{\varepsilon }\bigl(v_{\Gamma , \varepsilon }(s)+m_{0}\bigr)\bigr|_{H_{\Gamma }}^{2}ds\leq M_{6}. \label{es6} 
%\int_{0}^{T}\bigl|\boldsymbol{v}_{\varepsilon }(s)\bigr|_{\boldsymbol{W}}^{2}ds\leq M_{7}
\end{gather}
%{\it for a.a.\ $s\in (0, T)$. }

\paragraph{Proof}
From the definition of $\boldsymbol{\mu }_{\varepsilon }$, 
we can infer that 
\begin{eqnarray}
&&\mu _{\varepsilon }=\varepsilon \partial _{t}v_{\varepsilon }-\kappa _{1}\Delta v_{\varepsilon }
+\beta _{\varepsilon }(v_{\varepsilon }+m_{0})
-m(\boldsymbol{\beta }(\boldsymbol{v}_{\varepsilon }+m_{0}\boldsymbol{1}))+\widetilde{\pi }(v_{\varepsilon }+m_{0})-f \quad {\rm a.e.~in~}Q, \nonumber \\ \label{mu} \\
&&\mu _{\Gamma , \varepsilon }=\varepsilon \partial _{t}v_{\Gamma , \varepsilon }+\kappa _{1}\partial _{\boldsymbol{\nu }}v_{\varepsilon }
-\kappa _{2}\Delta _{\Gamma }v_{\Gamma , \varepsilon }
+\beta _{\Gamma , \varepsilon }(v_{\Gamma , \varepsilon }+m_{0})
-m(\boldsymbol{\beta }(\boldsymbol{v}_{\varepsilon }+m_{0}\boldsymbol{1}))\nonumber  \\
&&\quad \quad \quad +\widetilde{\pi }_{\Gamma }(v_{\Gamma , \varepsilon }+m_{0})-f_{\Gamma } \quad {\rm a.e.~on~}\Sigma . \label{mug}
\end{eqnarray}
Now, it follows from (\ref{es4}) that 
there exists a positive constant $\tilde{M}_{6}$ such that 
\begin{eqnarray}\label{volbeta}
&&\bigl|m\bigl(\boldsymbol{\beta }_{\varepsilon }\bigl(\boldsymbol{v}_{\varepsilon }(s)+m_{0}\boldsymbol{1}\bigr)\bigr)\bigr|^{2} \nonumber \\
&&\quad \leq \frac{2}{(|\Omega |+|\Gamma |)^{2}}\bigl(\bigl|\beta _{\varepsilon }\bigl(v_{\varepsilon }+m_{0}\bigr)\bigr|_{L^{1}(\Omega )}
+\bigl|\beta _{\Gamma , \varepsilon }\bigl(v_{\Gamma , \varepsilon }(s)+m_{0}\bigr)\bigr|_{L^{1}(\Gamma )}\bigr) \nonumber \\
&&\quad \leq \tilde{M}_{6}
\end{eqnarray}
for a.a.\ $s\in (0, T)$. 
Moreover, we test (\ref{mu}) at time $s\in (0, T)$ by $\beta _{\varepsilon }(v_\varepsilon (s)+m_{0})$ and exploit (\ref{mug}). 
Then, on account of the fact $(\beta _{\varepsilon }(v_{\varepsilon }+m_{0}))_{|_{\Gamma }}=\beta _{\varepsilon }(v_{\Gamma , \varepsilon }+m_{0})$, 
by integrating over $\Omega $ we deduce that 
\begin{eqnarray}
&&\kappa _{1}\int_{\Omega }\beta _{\varepsilon }'\bigl(v_{\varepsilon }(s)+m_{0}\bigr)\bigl|\nabla v_{\varepsilon }(s)\bigr|^{2}dx
+\kappa _{2}\int_{\Gamma }\beta _{\varepsilon }'\bigl(v_{\Gamma , \varepsilon }(s)+m_{0}\bigr)\bigl|\nabla _{\Gamma }v_{\Gamma , \varepsilon }(s)\bigr|^{2}d\Gamma \nonumber \\
&&\quad \quad \quad +\bigl|\beta _{\varepsilon }\bigl(v_{\varepsilon }(s)+m_{0}\bigr)\bigr|_{H}^{2}+\int_{\Gamma }\beta _{\Gamma , \varepsilon }\bigl(v_{\Gamma , \varepsilon }(s)+m_{0}\bigr)\beta _{\varepsilon }\bigl(v_{\Gamma , \varepsilon }(s)+m_{0}\bigr)d\Gamma  \nonumber \\
&&\quad \leq \bigl(f(s)+\mu _{\varepsilon }(s)-\varepsilon v_{\varepsilon }'(s)-\widetilde{\pi }\bigl(v_{\varepsilon }(s)+m_{0}\bigr), \beta _{\varepsilon }\bigl(v_{\varepsilon }(s)+m_{0}\bigr)\bigr)_{H} \nonumber \\
&&\quad \quad \quad +\bigl(m\bigl(\boldsymbol{\beta }_{\varepsilon }\bigl(\boldsymbol{v}_{\varepsilon }(s)+m_{0}\boldsymbol{1}\bigr)\bigr), \beta _{\varepsilon }(v_{\varepsilon }(s)+m_{0})\bigr)_{H} \nonumber \\
&&\quad \quad \quad \quad \quad +\bigl(f_{\Gamma }(s)+\mu _{\Gamma , \varepsilon }(s)-\varepsilon v_{\Gamma , \varepsilon }'(s)-\widetilde{\pi }_{\Gamma }\bigl(v_{\Gamma , \varepsilon }(s)+m_{0}\bigr), \beta _{\varepsilon }\bigl(v_{\Gamma , \varepsilon }(s)+m_{0}\bigr)\bigr)_{H_{\Gamma }} \nonumber \\
&&\quad \quad \quad \quad \quad \quad \quad +\bigl(m\bigl(\boldsymbol{\beta }_{\varepsilon }\bigl(\boldsymbol{v}_{\varepsilon }(s)+m_{0}\boldsymbol{1}\bigr)\bigr), \beta _{\varepsilon }(v_{\Gamma , \varepsilon }(s)+m_{0})\bigr)_{H_{\Gamma }} \label{intomg}
\end{eqnarray}
for a.a.\ $s\in (0, T)$. 
Now, from (\ref{betaine}), since the both sign of $\beta _{\varepsilon }(r)$ and $\beta _{\Gamma , \varepsilon }(r)$ is same for all $r\in \mathbb{R}$, 
we infer that 
\begin{eqnarray}
\int_{\Gamma }\beta _{\Gamma , \varepsilon }\bigl(v_{\Gamma , \varepsilon }(s)+m_{0}\bigr)\beta _{\varepsilon }\bigl(v_{\Gamma , \varepsilon }(s)+m_{0}\bigr)d\Gamma 
&=&\int_{\Gamma }\bigl|\beta _{\Gamma , \varepsilon }\bigl(v_{\Gamma , \varepsilon }(s)+m_{0}\bigr)\bigr|\bigl|\beta _{\varepsilon }\bigl(v_{\Gamma , \varepsilon }(s)+m_{0}\bigr)\bigr|d\Gamma \nonumber \\
&\geq &\frac{1}{2\rho }\int_{\Gamma }\bigl|\beta _{\varepsilon }\bigl(v_{\Gamma , \varepsilon }(s)+m_{0}\bigr)\bigr|^{2}d\Gamma 
-\frac{c_{0}^{2}}{2\rho }|\Gamma |. \label{sign}
\end{eqnarray}
Also, it holds 
\begin{equation}\label{intsei}
\int_{\Omega }\beta _{\varepsilon }'\bigl(v_{\varepsilon }(s)+m_{0}\bigr)\bigl|\nabla v_{\varepsilon }(s)\bigr|^{2}dx\geq 0, 
\quad \int_{\Gamma }\beta _{\Gamma , \varepsilon }'\bigl(v_{\Gamma , \varepsilon }(s)+m_{0}\bigr)\bigl|\nabla _{\Gamma }v_{\Gamma , \varepsilon }(s)\bigr|^{2}d\Gamma \geq 0. 
\end{equation}
Moreover, by using the {Y}oung inequality, the {L}ipschitz continuity of $\widetilde{\pi }, \widetilde{\pi }_{\Gamma }$
and (\ref{volbeta}), 
there exists a positive constant $\hat{M}_{6}$ such that 
\begin{eqnarray}
&&\bigl(f(s)+\mu _{\varepsilon }(s)-\varepsilon v_{\varepsilon }'(s)-\widetilde{\pi }\bigl(v_{\varepsilon }(s)+m_{0}\bigr), \beta _{\varepsilon }\bigl(v_{\varepsilon }(s)+m_{0}\bigr)\bigr)_{H} \nonumber \\
&&\quad \quad \quad +\bigl(m\bigl(\boldsymbol{\beta }_{\varepsilon }\bigl(\boldsymbol{v}_{\varepsilon }(s)+m_{0}\boldsymbol{1}\bigr)\bigr), \beta _{\varepsilon }\bigl(v_{\varepsilon }(s)+m_{0}\bigr)\bigr)_{H} \nonumber \\
&&\quad \leq \frac{1}{2}\bigl|\beta _{\varepsilon }\bigl(v_{\varepsilon }(s)+m_{0}\bigr)\bigr|_{H}^{2}
+4\bigl|f(s)\bigr|_{H}^{2}+4\bigl|\mu _{\varepsilon }(s)\bigr|_{H}^{2}
+4\varepsilon ^{2}\bigl|v_{\varepsilon }'(s)\bigr|_{H}^{2}
+4\bigl|\widetilde{\pi }\bigl(v_{\varepsilon }(s)+m_{0}\bigr)\bigr|_{H}^{2} \nonumber \\
&&\quad \quad \quad +\bigl|m\bigl(\boldsymbol{\beta }_{\varepsilon }\bigl(\boldsymbol{v}_{\varepsilon }(s)+m_{0}\boldsymbol{1}\bigr)\bigr)\bigr|_{H}^{2} \nonumber \\
&&\quad \leq \frac{1}{2}\bigl|\beta _{\varepsilon }\bigl(v_{\varepsilon }(s)+m_{0}\bigr)\bigr|_{H}^{2}
+\hat{M}_{6}\bigl(\bigl|f(s)\bigr|_{H}^{2}+\bigl|\mu _{\varepsilon }(s)\bigr|_{H}^{2}
+\varepsilon ^{2}\bigl|v_{\varepsilon }'(s)\bigr|_{H}^{2}+\bigl|v_{\varepsilon }(s)\bigr|_{H}^{2}+1\bigr) \nonumber \\
&&\quad \quad \quad +|\Omega |\tilde{M}_{6} \label{homg}
\end{eqnarray}
and 
\begin{eqnarray}
&&\bigl(f_{\Gamma }(s)+\mu _{\Gamma , \varepsilon }(s)-\varepsilon v_{\Gamma , \varepsilon }'(s)-\widetilde{\pi }_{\Gamma }\bigl(v_{\Gamma , \varepsilon }(s)+m_{0}\bigr), \beta _{\varepsilon }\bigl(v_{\Gamma , \varepsilon }(s)+m_{0}\bigr)\bigr)_{H_{\Gamma }} \nonumber \\
&&\quad \quad \quad +\bigl(m\bigl(\boldsymbol{\beta }_{\varepsilon }\bigl(\boldsymbol{v}_{\varepsilon }(s)+m_{0}\boldsymbol{1}\bigr)\bigr), \beta _{\varepsilon }(v_{\Gamma , \varepsilon }(s)+m_{0})\bigr)_{H} \nonumber \\
&&\quad \leq \frac{1}{4\rho }\bigl|\beta _{\varepsilon }\bigl(v_{\Gamma , \varepsilon }(s)+m_{0}\bigr)\bigr|_{H_{\Gamma }}^{2}
+2\rho \bigl|f_{\Gamma }(s)\bigr|_{H}^{2}+2\rho \bigl|\mu _{\Gamma , \varepsilon }(s)\bigr|_{H_{\Gamma }}^{2}
+2\rho \varepsilon ^{2}\bigl|v_{\Gamma , \varepsilon }'(s)\bigr|_{H_{\Gamma }}^{2} \nonumber \\
&&\quad \quad \quad +2\rho \bigl|\widetilde{\pi }_{\Gamma }\bigl(v_{\Gamma , \varepsilon }(s)+m_{0}\bigr)\bigr|_{H_{\Gamma }}^{2}
+2\rho \bigl|m\bigl(\boldsymbol{\beta }_{\varepsilon }\bigl(\boldsymbol{v}_{\varepsilon }(s)+m_{0}\boldsymbol{1}\bigr)\bigr)\bigr|_{H_{\Gamma }}^{2} \nonumber \\
&&\quad \leq \frac{1}{4\rho }\bigl|\beta _{\varepsilon }\bigl(v_{\Gamma , \varepsilon }(s)+m_{0}\bigr)\bigr|_{H_{\Gamma }}^{2} 
+2\rho |\Gamma |\tilde{M}_{6} \nonumber \\
&&\quad \quad \quad +\rho \hat{M}_{6}\bigl(\bigl|f_{\Gamma }(s)\bigr|_{H_{\Gamma }}^{2}+\bigl|\mu _{\Gamma , \varepsilon }(s)\bigr|_{H_{\Gamma }}^{2}
+\varepsilon ^{2}\bigl|v_{\Gamma , \varepsilon }'(s)\bigr|_{H_{\Gamma }}^{2}
+\bigl|v_{\Gamma , \varepsilon }(s)\bigr|_{H_{\Gamma }}^{2}+1\bigr) \label{hgomg}
\end{eqnarray}
for a.a.\ $s\in (0, T)$. 
Thus, from Lemma 3.3, 3.4 and (\ref{PW0}), by combining from (\ref{intomg})--(\ref{hgomg}) and integrating them over $(0, T)$, we can conclude existence of the constant $M_{6}$ satisfying (\ref{es6}). 
\hfill $\Box $

%and we already obtained estimates respect to $\boldsymbol{v}_{\varepsilon }$ and $\boldsymbol{\beta }_{\varepsilon }$ enough to prove. 
%The key point to prove lemmas is that we can adapt the theory of the elliptic regularity (see, e.g.\, \cite[Theorem 3.2, p.\ 1.79]{BG87})
%becuse we have the term $\Delta _{\Gamma }v _{\varepsilon }$ in the approximate sequence $\mu _{\Gamma , \varepsilon }$. 
%This is the big advantage of dynamic boundary condition (\ref{eq5}). 
%Therefore, we omit the proof. 

%%%seventh estimate%%%%%%%%%%%%%%%%%%%%%%%%%%%%%%%%%%%%%%%%%%
\paragraph{Lemma 3.8.}
{\it There exists a positive constant $M_{7}$, independent of 
$\varepsilon \in (0, 1]$, such that }
\begin{equation}\label{es7}
\kappa _{1}\int_{0}^{T}\bigl|\Delta v_{\varepsilon }(s)\bigr|_{H}^{2}ds
+\int_{0}^{T}\bigl|v_{\varepsilon }(s)\bigr|_{H^{\frac{3}{2}}(\Omega )}^{2}ds
+\int_{0}^{T}\bigl|\partial _{\boldsymbol{\nu }}v_{\varepsilon }(s)\bigr|_{H_{\Gamma }}^{2}ds
\leq M_{7}. 
\end{equation}
%{\it for a.a.\ $s\in (0, T)$. }
%\paragraph{Proof}
This lemma is proved exactly the same as in \cite[Lemmas 4.4]{CF15}
because the necessary uniform estimates to prove it is obtained by Lemmas 3.3, 3.4, 3.6 and 3.7. 
Sketching simply, comparing in (\ref{mu}) we deduce that $|\Delta v_{\varepsilon }|_{L^{2}(0, T; H)}$ is uniformly bounded. 
Moreover, by using the theory of the elliptic regularity and the trace theory
(see e.g., \cite[Theorem 3.2, p.\ 1.79 and Theorem 2.25, p.\ 1.62]{BG87}, respectively), 
we can conclude that (\ref{es7}) holds.

%%%eighth estimate%%%%%%%%%%%%%%%%%%%%%%%%%%%%%%%%%%%%%%%%%%
\paragraph{Lemma 3.9.}
{\it There exists a positive constant $M_{8}$, independent of 
$\varepsilon \in (0, 1]$, such that }
\begin{equation}\label{es8}
\int_{0}^{T}\bigl|\beta _{\Gamma , \varepsilon }\bigl(v_{\Gamma , \varepsilon }(s)+m_{0}\bigr)\bigr|_{H_{\Gamma }}^{2}ds
\leq M_{8}.
\end{equation}
%{\it for a.a.\ $s\in (0, T)$. }

\paragraph{Proof}
We test (\ref{mug}) at time $s\in (0, T)$ by $\beta _{\Gamma , \varepsilon }(v_{\Gamma , \varepsilon }(s)+m_{0})$ and integrating it over $\Gamma $. 
Then, by using the {Y}oung inequality and the {L}ipschitz continuity of $\widetilde{\pi }_{\Gamma }$, 
there exists a positive constant $\tilde{M}_{8}$ such that 
\begin{eqnarray}
&&\kappa _{2}\int_{\Gamma }\beta _{\Gamma , \varepsilon }'\bigl(v_{\Gamma , \varepsilon }(s)+m_{0}\bigr)\bigl|\nabla _{\Gamma }v_{\Gamma , \varepsilon }(s)\bigr|^{2}d\Gamma 
+\bigl|\beta _{\Gamma , \varepsilon }\bigl(v_{\Gamma , \varepsilon }(s)+m_{0}\bigr)\bigr|_{H_{\Gamma }}^{2} \nonumber \\
&&\quad =\bigl(f_{\Gamma }(s)+\mu _{\Gamma }(s)-\varepsilon v_{\Gamma , \varepsilon }'(s)-\partial _{\boldsymbol{\nu }}v_{\varepsilon }(s)
-\widetilde{\pi }_{\Gamma }\bigl(v_{\Gamma , \varepsilon }(s)+m_{0}\bigr), \beta _{\Gamma , \varepsilon }\bigl(v_{\Gamma , \varepsilon }(s)+m_{0}\bigr)\bigr)_{H_{\Gamma }} \nonumber \\
&&\quad \leq \frac{1}{2}\bigl|\beta _{\Gamma , \varepsilon }\bigl(v_{\Gamma , \varepsilon }(s)+m_{0}\bigr)\bigr|_{H_{\Gamma }}^{2} 
+\bigl|m\bigl(\boldsymbol{\beta }_{\varepsilon }\bigl(\boldsymbol{v}_{\varepsilon }(s)+m_{0}\bigr)\bigr)\bigr|_{H_{\Gamma }}^{2} \nonumber \\
&&\quad \quad \quad +\tilde{M}_{8}\bigl(\bigl|f_{\Gamma }(s)\bigr|_{H_{\Gamma }}^{2}
+\bigl|\mu _{\Gamma }(s)\bigr|_{H_{\Gamma }}^{2}
+\varepsilon ^{2}\bigl|v_{\Gamma , \varepsilon }'(s)\bigr|_{H_{\Gamma }}^{2}
+\bigl|\partial _{\boldsymbol{\nu }}v_{\varepsilon }(s)\bigr|_{H_{\Gamma }}^{2}
+\bigl|v_{\Gamma , \varepsilon }(s)\bigr|_{H_{\Gamma }}^{2}+1\bigr) \nonumber \\
&&\quad \leq \frac{1}{2}\bigl|\beta _{\Gamma , \varepsilon }\bigl(v_{\Gamma , \varepsilon }(s)+m_{0}\bigr)\bigr|_{H_{\Gamma }}^{2} 
+|\Gamma |\tilde{M}_{6} \nonumber \\
&&\quad \quad \quad +\tilde{M}_{8}\bigl(\bigl|f_{\Gamma }(s)\bigr|_{H_{\Gamma }}^{2}
+\bigl|\mu _{\Gamma }(s)\bigr|_{H_{\Gamma }}^{2}
+\varepsilon ^{2}\bigl|v_{\Gamma , \varepsilon }'(s)\bigr|_{H_{\Gamma }}^{2}
+\bigl|\partial _{\boldsymbol{\nu }}v_{\varepsilon }(s)\bigr|_{H_{\Gamma }}^{2}
+\bigl|v_{\Gamma , \varepsilon }(s)\bigr|_{H_{\Gamma }}^{2}+1\bigr) \nonumber \\ \label{es8ine}
\end{eqnarray}
for a.a.\ $s\in (0, T)$. 
Noting that it holds 
\begin{equation*}
\kappa _{2}\int_{\Gamma }\beta _{\Gamma , \varepsilon }'\bigl(v_{\Gamma , \varepsilon }(s)+m_{0}\bigr)\bigl|\nabla _{\Gamma }v_{\Gamma , \varepsilon }(s)\bigr|^{2}d\Gamma \geq 0. 
\end{equation*}
Thus, on account of Lemma 3.3, 3.4, 3.6 and 3.8, 
by integrating (\ref{es8ine}) over $(0, T)$, 
we can find a positive constant $M_{7}$ such that the estimate (\ref{es8}) holds. 
\hfill $\Box $

%%%nineth estimate%%%%%%%%%%%%%%%%%%%%%%%%%%%%%%%%%%%%%%%%%%
\paragraph{Lemma 3.10.}
{\it There exists a positive constant $M_{9}$, independent of 
$\varepsilon \in (0, 1]$, such that }
\begin{equation*}\label{es9}
\int_{0}^{T}\bigl|\boldsymbol{v}_{\varepsilon }(s)\bigr|_{\boldsymbol{W}}^{2}ds
\leq M_{9}. 
\end{equation*}
%{\it for a.a.\ $s\in (0, T)$. }
%\paragraph{Proof}
This lemma is also proved the same as in \cite[Lemmas 4.5]{CF15}. 
The key point to prove it is that we can obtain the uniform estimates of 
$|\Delta _{\Gamma }v_{\Gamma , \varepsilon }|_{L^{2}(0, T; H_{\Gamma })}$ 
by comparing in (\ref{mug}). 
We omit the proof.

%\int_{0}^{T}\bigl|\boldsymbol{\mu }_{\varepsilon }(s)\bigr|_{\boldsymbol{V}}^{2}ds\leq M_{5}, \\
%\int_{0}^{T}\bigl|\boldsymbol{\beta }_{\varepsilon }\bigl(u_{\varepsilon }(s)\bigr)\bigr|_{H}^{2}ds
%+\int_{0}^{T}\bigl|\boldsymbol{\beta }_{\Gamma , \varepsilon }\bigl(u_{\Gamma , \varepsilon }(s)\bigr)\bigr|_{H_{\Gamma }}^{2}ds\leq M_{6}, \\
%\int_{0}^{T}\bigl|\boldsymbol{v}_{\varepsilon }(s)\bigr|_{\boldsymbol{W}}^{2}ds\leq M_{7}

%%%%% Section 4. %%%%%

\section{Proof of convergence theorem}
%\subsection{Passage to the limit }
\indent

In this section, 
we obtain the existence of weak solution of (P) by performing passage to the limit for the approximate problem (P)$_{\varepsilon }$. 
The convergence theorem is also nearly the same \cite[Sect.\ 4]{CF15}. 
The different point from \cite{CF15} is that the component of the weak solution of (P) satisfies 
(\ref{weak2}) and the periodic property (\ref{ic}).

Thanks to the previous estimates Lemmas from 3.3 to 3.10, 
there exist a subsequence $\{\varepsilon _{k}\}_{k\in \mathbb{N}}$ with $\varepsilon _{k}\to 0$ as $k\to \infty $ and some limits functions 
$\boldsymbol{v}\in H^{1}(0, T; \boldsymbol{V}_{0}^{*})\cap L^{\infty }(0, T; \boldsymbol{V}_{0})\cap L^{2}(0, T; \boldsymbol{W})$, 
$\boldsymbol{\mu }\in H^{1}(0, T; \boldsymbol{V})$, 
$\xi \in L^{2}(0, T; H)$
and $\xi _{\Gamma }\in L^{2}(0, T; H_{\Gamma })$ such that 
\begin{equation}\label{limv}
\boldsymbol{v}_{\varepsilon _{k}}\to \boldsymbol{v} \quad {\rm weakly~star~in~} H^{1}(0, T; \boldsymbol{V}_{0}^{*})\cap L^{\infty }(0, T; \boldsymbol{V}_{0})\cap L^{2}(0, T; \boldsymbol{W}), 
\end{equation}
\begin{equation*}\label{limepv}
\varepsilon _{k}\boldsymbol{v}_{\varepsilon _{k}}\to 0 \quad {\rm strongly~in~} H^{1}(0, T; \boldsymbol{H}_{0}), 
\end{equation*}
\begin{equation*}\label{limmu}
\boldsymbol{\mu }_{\varepsilon _{k}}\to \boldsymbol{\mu } \quad {\rm weakly~in~} L^{2}(0, T; \boldsymbol{V}), 
\end{equation*}
\begin{equation}\label{limbe}
\beta _{\varepsilon _{k}}(u_{\varepsilon _{k}})\to \xi \quad {\rm weakly~in~} L^{2}(0, T; H), 
\end{equation}
\begin{equation}\label{limbeg}
\beta _{\Gamma , \varepsilon _{k}}(u_{\Gamma , \varepsilon _{k}})\to \xi _{\Gamma } \quad {\rm weakly~in~} L^{2}(0, T; H_{\Gamma })
\end{equation}
as $k\to \infty $. 
Owing to (\ref{limv}) and a well-known compactness results (see e.g., \cite{Sim87}), we obtain 
\begin{equation}\label{strov}
\boldsymbol{v}_{\varepsilon _{k}}\to \boldsymbol{v} \quad {\rm strongly~in~} C([0, T]; \boldsymbol{H}_{0})\cap L^{2}(0, T; \boldsymbol{V}_{0})
\end{equation}
as $k\to \infty $. 
This yeilds that 
\begin{equation}\label{strou}
\boldsymbol{u}_{\varepsilon _{k}}\to \boldsymbol{u}:=\boldsymbol{v}+m_{0}\boldsymbol{1} \quad {\rm strongly~in~} C([0, T]; \boldsymbol{H}_{0})\cap L^{2}(0, T; \boldsymbol{V}_{0})
\end{equation}
as $k\to \infty $. 
Therefore, from (\ref{strou}) and the {L}ipschitz continuity of $\widetilde{\pi }, \widetilde{\pi }_{\Gamma }$, we deduce that 
\begin{equation*}
\widetilde{\boldsymbol{\pi }}(\boldsymbol{u}_{\varepsilon _{k}})\to \widetilde{\boldsymbol{\pi }}(\boldsymbol{u}) \quad {\rm strongly~in~} C([0, T]; \boldsymbol{H}). 
\end{equation*}
Hence, by passing to the limit in (\ref{apromu}) and (\ref{vmu}), 
we obtain (\ref{weak1}) and the following weak formulation: 
\begin{gather}
%\bigl\langle \boldsymbol{v}'(t), \boldsymbol{z}\bigr\rangle _{\boldsymbol{V}_{0}^{*}, \boldsymbol{V}_{0}}+a\bigl(\boldsymbol{\mu }(t), \boldsymbol{z}\bigr)=0 \quad {\rm for~all~}z\in \boldsymbol{V}_{0} \nonumber \\
\bigl(\boldsymbol{\mu }(t), \boldsymbol{z}\bigr)_{\boldsymbol{H}}=a\bigl(\boldsymbol{v}(t), \boldsymbol{z}\bigr)+\bigl(\boldsymbol{\xi }(t)-m\bigl(\boldsymbol{\xi }(t)\bigr)\boldsymbol{1}+\widetilde{\boldsymbol{\pi }}\bigl(\boldsymbol{v}(t)\bigr)-\boldsymbol{f}, \boldsymbol{z}\bigr)_{\boldsymbol{H}} \quad {\rm for~all~}z\in \boldsymbol{V} \label{weaktil}
\end{gather}
for a.a.\ $t\in (0, T)$ where $\boldsymbol{\xi }:=(\xi , \xi _{\Gamma })$, 
because of the property (\ref{projection}) of linear bounded operator $\boldsymbol{P}$. 
Now, we can infer 
$v+m_{0}\in D(\beta )$ and $v_{\Gamma }+m_{0}\in D(\beta _{\Gamma })$. 
Hence, from the form (\ref{pi}) and (\ref{pig}), we deduce that 
%$\widetilde{\boldsymbol{\pi }}(\boldsymbol{v}+m_{0}\boldsymbol{1})=\boldsymbol{\pi }(\boldsymbol{v}+m_{0}\boldsymbol{1})$, 
%that is, 
$\widetilde{\pi }(v+m_{0})=\pi (v+m_{0})$  a.e.\ in $Q$ and 
$\widetilde{\pi }_{\Gamma }(v+m_{0})=\pi _{\Gamma }(v_{\Gamma }+m_{0})$ a.e.\ on $\Sigma $. 
This implies that we obtain (\ref{weak2}) replaced by (\ref{weaktil}). 
Moreover, it follows from (\ref{strov}) that 
$$\boldsymbol{v}(0)=\boldsymbol{v}(T) \quad {\rm in~}\boldsymbol{H}_{0}. $$
Also, due to (\ref{limbe}), (\ref{limbeg}), (\ref{strou}) and the monotonicity of $\beta $, 
from the fact \cite[Prop.\ 2.2, p.\ 38]{Bar10} we obtain 
\begin{equation*}
\xi \in \beta (v+m_{0}) \quad {\rm a.e.~in~}Q, \quad \xi _{\Gamma }\in \beta _{\Gamma }(v+m_{0}) \quad {\rm a.e.~on~}\Sigma . 
\end{equation*}
Thus, we complete the proof of Theorem 2.1. 
%\hfill  $\Box $

%%%%% Section 5. %%%%%
%\section{Choosing prototype double well potential}
%\setcounter{equation}{0}
%\indent

\paragraph{Remark 4.1.}
In the previous sections, we impose the restricted assumption (A5) respect to the domains $D(\beta )$ and $D(\beta _{\Gamma })$. 
This is an essential assumption if we treat the multivalued $\beta , \beta _{\Gamma }$. 
However, refferring to \cite{LWZ16}, we can avoid the assumption (A5) when focusing only on prototype double potential, 
namely when we choose $\beta (u)=u^{3}$, $\beta _{\Gamma }(u_{\Gamma })=u_{\Gamma }^{3}$, 
$\pi (u)=-u$ and $\pi _{\Gamma }(u_{\Gamma })=-u_{\Gamma }$ in (\ref{eq2}) and (\ref{eq5}). 
The reason why we can avoid it is that we need not to use (\ref{pibdd}) to obtain estimates. 
The method to do so is used in \cite{LWZ16}. 
Concretely, it is using the H\"older inequality mainly. 
We can employ the same method even if we consider the problem 
with dynamic boundary condition. 
Therefore, we can prove Theorem 2.1 without the assumption (A5) in the case of the prototype double potential.

\end{document}